\documentclass[12pt,a4paper]{article}

\usepackage{amsmath,amsthm,amsfonts,amssymb,amscd}
\usepackage{graphicx}

\newtheorem{theorem}{Theorem}[section]
\newtheorem{proposition}{Proposition}[section]
\newtheorem{lemma}{Lemma}[section]

\newcommand{\erm}{{\rm e}}

\newcommand{\cinfty} {\hbox{\rm C}^\infty}

\newcommand{\ov} {\overline}

\newcommand{\Real}{{\rm Re}}

\newcommand{\Oc} {{\cal O}}

\newcommand{\calG}{{\cal G}}

\newcommand{\dt} {\delta}
\newcommand{\Dt} {\Delta}
\newcommand{\ep} {\epsilon}
\newcommand{\Gm} {\Gamma}

\newcommand{\lb} {\lambda}

\newcommand{\sg} {\sigma}

\def\S{{\mathbb S}}
\def\Z{{\mathbb Z}}

\def\R{{\mathbb R}}

\def\C{{\mathbb C}}

\newcommand{\bproof}{\noindent {\it Proof. }}
\newcommand{\eproof}{\hfill$\Box$}
\newcommand{\nd} {\noindent}

\begin{document}
\begin{titlepage}

\begin{center}

  {\large \bf Vortex   on surfaces and 
    Brownian-motion 
  in higher dimensions: special metrics}

\end{center}

\vskip  1.0truecm
\centerline {{\large Clodoaldo Grotta-Ragazzo$^*$.}}

\vskip 0.5truecm



\begin{abstract} 
  A single hydrodynamic vortex on a  surface  will in general  moves unless its 
  Riemannian metric is a  special
  ``Steady Vortex Metric'' (SVM).  Metrics 
  of constant curvature are SVM only in surfaces of genus zero  and one.

In this paper:

\begin{enumerate}
  \item I show that  K. Okikiolu's work on the regularization of the spectral zeta function leads to the conclusion that each conformal class of every compact surface with a genus of two or more possesses at least one Steady Vortex Metric (SVM).

  \item I apply a probabilistic interpretation of the regularized zeta function for surfaces, as developed by P. G. Doyle and J. Steiner, to extend the concept of SVM to higher dimensions.
\end{enumerate}

 The new special metric, which aligns with the Steady Vortex Metric (SVM) in two dimensions, has been termed the ``Uniform Drainage Metric'' for the following reason: For a compact Riemannian manifold \( M \), the ``narrow escape time'' (NET) is defined as the expected time for a Brownian motion starting at a point \( p \) in \( M \setminus B_\epsilon(q) \) to remain within this region before escaping through the small ball \( B_\epsilon(q) \), which is centered at \( q \) with radius \( \epsilon \) and acts as the escape window. The manifold is said to possess a uniform drainage metric if, and only if, the spatial average of NET, calculated across a uniformly distributed set of initial points \( p \), remains invariant regardless of the position of the escape window \( B_\epsilon(q) \), as \( \epsilon \) approaches \( 0 \).

\end{abstract}

\vskip .5truecm

\nd
{\bf Key words: point vortex, Riemann surfaces, diffusion process, Brownian motion,  special metrics,
spectral zeta function.} 
\vskip .5truecm

\nd
{\bf AMS Classification: 76B47, 30F30, 58J65, 31C12, 60J45, 53C25}

\vskip .5truecm
\nd{\bf Abbreviated title: SVM and  Brownian-motion 
  special metrics}

\vfill
\hrule
\smallskip
\nd
$^*$ Instituto de Matem\'{a}tica e Estat\'{i}stica da
Universidade de S\~{a}o Paulo,\\
\nd Rua do Mat\~ao 1010, 05508-090, S\~{a}o Paulo, SP, Brazil.\\
\nd Partially supported by FAPESP grant 2016/25053-8.\\
 \nd email: ragazzo@usp.br \\
\nd ORCID: 0000-0002-4277-4173
\end{titlepage}

\pagebreak
\section{Introduction}

\label{eulersec}

The motion of point vortices on the  plane  is a classical
subject in fluid mechanics that goes back to Helmholtz, Kelvin, and
Kirchhoff.  The first to consider the motion of
point vortices on a curved surface, the sphere embedded in $\R^3$,
was Zermello in 1902.
The paper \cite{borigid}    
has a historical review 
on the early research on hydrodynamic vortices on surfaces. 
An intrinsic definition of the motion of vortices on a surface,
which is independent of the embedding
of the surface in  $\R^3$ and on coordinates, started with Boatto and Koiller   
\cite{boatto2008vortices} (see also \cite{kb2}
\cite{db}  \cite{ragvil}) and was recently completed 
by Bj\"orn Gustafsson \cite{gustafsson2019vortex} \cite{gustafsson2022dipole}.

A single vortex in the Euclidean plane, or in the round sphere, or in a flat torus
does not move, and this motivated  the definition of ``Steady Vortex Metric'' 
 \cite{ragvil}: a Riemannian metric for
 which a single vortex does not move regardless of its position.
J. Koiller conjectured that  a
single vortex in a compact surface of constant curvature
and of a genus greater than one does move. In  \cite{ragazzo2017motion},
\cite{errata}
Koiller's conjecture
was numerically verified for a particular surface of
constant curvature of genus two:
the Bolza surface.
 This result  motivated the first main
question to be answered in this work:
Does a steady vortex metric
exist  on any orientable compact surface of a genus
greater than one?

K. Okikiolu proved that a certain functional on the space of Riemannian metrics,
which is an analog  for closed
 surfaces of the ADM mass from general relativity, has a minimizing
 metric on each conformal class.
 It turns out that the special metrics of  Okikiolu are  steady vortex metrics,
 which gives a positive answer to the question in the paragraph
 above. 
This raises the question about the ``meaning'' (or  properties)
of this  special metric. The steady vortex  metric minimizes a certain functional
\cite{okigenus} and has the property
 in its name, but does it have
any other interesting geometrical property besides those? This question
was the second motivation for this work. 

The special metric found by Okikiolu is a critical point of a functional
related to the regularized Green's function of the Laplacian: the ``Robin
function''.  P. G. Doyle and J. Steiner \cite{doyle2017spectral}
gave a probabilistic interpretation to the Robin function that is
related to the concept of
``Narrow-Escape-Time''(NET) \cite{holcman2014narrow}.
The NET
  is defined as the expected time for a Brownian motion starting at
   $p$  in 
 $ M\backslash B_\epsilon(q)$
 to remain within this region before escaping through the small ball $B_\epsilon(q)$,
 which is centered at 
 $q$ with radius 
 $\epsilon$ and acts as the escape window.

 The NET is an  important abstraction  in science, as
 argued by Holcman and Schuss in the Introduction of \cite{holcman2014narrow}:
   ``The  narrow escape problem in diffusion theory, which goes back
to Helmholtz (Helmholtz (1860)) and Lord Rayleigh (Rayleigh (1945)) in the context
of the theory of sound, is to calculate the mean first passage time
of Brownian motion to a small absorbing window....
The renewed interest in the problem
is due to the emergence of the narrow escape time (NET) as a key to the
determination of biological cell function from its geometrical structure.
The NET is ubiquitous in molecular and cellular
biology and is manifested in stochastic models of chemical reactions...''

The average NET, with respect to a uniform distribution 
of  initial positions (volume measure),
that a particle takes to escape from 
$S\backslash B_\epsilon(q)$ through the small window $B_\epsilon(q)$ is proportional
to $-\log \epsilon+ R(q)+\Oc(\epsilon)$, where $R$ is the Robin function.
So, for small $\epsilon$ the Robin function
indicates the drainage capacity  of different points $q$ in $S$.
 The Robin function is constant if, and only if, the metric is a
 steady vortex metric (SVM). Therefore, in a surface with a  $SVM$ the drainage
 capacity of different points is the same and this lead to the alternative
 name ``{\it uniform drainage metric}'', a property that makes sense in dimensions
 larger than two. Note: the notion of hydrodynamic {\it point} vortex cannot be
 generalized to dimensions
greater than two.

The main contribution in this paper is the definition of
uniform drainage metric in  dimensions greater than two and
its geometric characterization in dimensions 3 and 4.

 Following the same steps given in this paper, a characterization of a uniform drainage metric in higher dimensions can be accomplished by means of certain coefficients that appear in the so-called Minakshisundaram-Plejel asymptotic expansion of the heat kernel. I prefer not to state any results in this direction because, in higher dimensions, it is necessary to compute more of these coefficients, which can be expressed in terms of powers of the Laplacian and the distance function $\ell$, and they become very complicated \cite{polt}.

The existence of uniform drainage surfaces of arbitrary finite genus in
any conformal class 
is guaranteed by the theorem of Okikiolu. In higher dimensions any compact 
Riemannian manifold that is a homogeneous space is a uniform drainage
manifold.\footnote{ There is a  special class of Riemannian metrics on
  closed manifolds
  that are critical metrics of
  the trace of the heat kernel under conformal variations of the metric
  \cite{el2002critical}. A metric in this special 
  class is always a uniform drainage metric (a consequence of Theorem 4.1 (ii)
  in  \cite{el2002critical}). The metric
  of any Riemannian homogeneous space is critical for the trace of the heat
  kernel.}  
Does there exist a closed
(compact and boundaryless) manifold that does not admit a uniform drainage metric?

 This paper is organized as follows.

 In Section \ref{existence} I give a precise definition
of  the steady vortex metric  and present  two fundamental theorems
that stem
from Okikiolu's work. I then use these theorems to compare the
steady vortex metric with other natural Riemannian  metrics:
of constant curvature, canonical or Bergman, and Arakelov.
The proofs  of the two theorems are presented in  Appendix
\ref{proofsec} in a slightly different way than those given by Okikiolu.
These theorems  plus some simple arguments imply:
``No orientable surface of genus 2 and of constant curvature
is a Steady Vortex Surface.''

In Section \ref{hdim}, I present a regularization of the Green's function
in dimensions greater than two using
the Minakshisundaram-Plejel asymptotic expansion of the   heat kernel.
This provides a
definition of the Robin function in  higher dimensions.
In the Appendix \ref{min} I show that
the Robin function can be written in terms of the analytic
extension of the Minakshisundaram–Pleijel zeta function, and
therefore  uniform drainage manifolds have a special spectral
property derived from this relation. The relation between the Robin
function and the  Minakshisundaram–Pleijel zeta function 
appeared in  \cite{steiner2005geometrical}, for surfaces,
and in  \cite{bilal2013multi}, in a more general context and in dimension
greater than  two.

In Section \ref{netsec} I give a characterization
of uniform drainage metric in dimensions
2, 3 and 4. In dimension 4 a uniform drainage metric has constant
Robin function (a global property) and constant  scalar curvature
(a local property).

In Section \ref{torsec}, I present a family of non-flat tori
found by Okikiolu \cite{okitorus},
and which will be called Okikiolu's tori,
that are
uniform drainage surfaces. These tori are the only non-constant
curvature uniform drainage surface that are  explicitly known.
For any $a>\sqrt{\pi/2}$ there is an Okikiolu's torus that is
conformally equivalent to the flat  torus
$\R^2/(a\Z\times a^{-1}\Z)$.
Therefore,  uniform drainage metrics may not be unique  in a
conformal class. The curvature of the Okikiolu's  tori
was computed in ibid., where it was realized that in the limit as $a\to\infty$
the curvature at  almost every  point of the  torus
tends to $1/\sqrt{4\pi}$. In Section  \ref{torsec}
I embed a cylinder in $\R^3$ whose
quotient under the group of translations along the cylinder axis is an
Okikiolu's torus. In this way, I can visualize the deformation
of a flat torus into a pinched torus that is isometric to a round sphere
with two opposite points being identified.  The deformation is done
along an interesting family  
of uniform drainage surfaces.

I finally remark about a possible upshot  of the  relation between the
Robin function and the drainage capacity of different points. The importance of the
NET in cellular biology is partially due to diffusion processes that occur in membranes
towards special
exit gates (escape windows). The minimum of the Robin function is an equilibrium position
of a single vortex \cite{errata} and also a point where the drainage capacity of the surface,
as defined above, is maximum. Equilibrium positions of systems of point vortices, an issue that
has been extensively studied, also have a probabilistic interpretation. If the position
of a vortex is related to an entrance or exit gate, depending on the vortex sign,
then some equilibrium 
configurations  will certainly be more efficient in connecting different gates
by means of  diffusion than others. If this idea is correct, then the 
the importance of equilibrium configurations on surfaces of spheres, including those
which are not round, will be greatly enhanced.

\section{Steady vortex metrics on
  orientable 
  closed surfaces.}
 \label{existence}

 The definition of a hydrodynamic vortex requires some preliminaries (see \cite{ragvil}). The fundamental equations of hydrodynamics on $S$, Euler's equations, necessitate that $S$ be endowed with a Riemannian metric $g$. Here, $g$ represents a smooth family of inner products on the tangent spaces of $S$. In local coordinates, the Riemannian metric is given by $g=\sum_{jk}g_{jk}dx_j\otimes dx_k$. The associated volume form is $\mu=\sqrt{|g|}dx_1\wedge dx_2$, where $|g|$ denotes the absolute value of the determinant of the matrix $g_{jk}$

 In a neighborhood of each point of \( S \), there exist coordinates (sometimes called isothermal coordinates) in which \( g=\lambda^2(x)(dx_1^2+dx_2^2) \) and \( \mu=\lambda^2(x)dx_1\wedge dx_2 \).
 The existence of isothermal coordinates is a manifestation of the fact that any surface is locally conformal to the Euclidean plane. In this paper, I will also use \( \lambda^2 \) to denote the conformal factor between arbitrary given metrics \( g_0 \) and \( g_1 \). This will be explicitly stated when used.

 The one-forms \( \theta_1=\lambda dx_1 \) and \( \theta_2=\lambda dx_2 \) constitute an orthonormal moving coframe. The Hodge-star operator acts linearly on forms and is defined by
\[
 \ast 1=\theta_1\wedge\theta_2=\mu,\qquad 
\ast\theta_1=\theta_2,\qquad \ast\theta_2=-\theta_1, \qquad
 \ast\mu=1.
\]
The Laplace operator acting on
functions is given by \( \Delta=\ast d\ast d=\frac{1}{\lambda^2}(\partial_{x_1}^2+\partial_{x_2}^2) \)
and the Gaussian curvature by \( K=-\frac{1}{\lambda^2}\Delta \log \lambda \).

Let $V=\int_S\mu$ be  the total area of $S$.
The  Green's function of $(S,g)$ is the unique  solution 
in distribution sense
to the equation 
\begin{equation}
  -\Dt_q G(q,p)=\dt_p(q)-V^{-1},
\label{compacteq}
\end{equation}
that has the following properties (see \cite{aubin}, theorem 4.13):
\begin{itemize}
  \item for all functions $\phi\in C^2$
  \begin{equation}
  \phi(p)=\frac{1}{V}\int_S\phi\mu - \int_SG(q,p)\Dt \phi(q)\mu(q)\label{G},
  \end{equation}
  \item $G(q,p)$ is $\cinfty$ on $S\times S$  minus the diagonal,
\item $G$ is symmetric $G(q,p)=G(p,q)$,
  \item $G$ is bounded from below
  and  $\int_S G(q,p)\mu(q)=0$.
\end{itemize}
A  point vortex of intensity $\Gm\in\R$ 
at the  point $p$   
is the fluid velocity field defined on $S-\{p\}$ given by
$q\to\ast\nabla \Gamma G(q,p)$, where $\nabla$ is the gradient
operator and $\ast$ is the operator that
rotates a vector by $\pi/2$.

The Robin function (the regularization of $G$)
is a  $\cinfty$ function on $S$ (\cite{ragvil} Theorem 5.1) defined as
\begin{equation}
R(p)=
\lim_{\ell(q,p)\to 0} \left[
G(q,p)+\frac{1}{2\pi}\log \ell(q,p)\right],
\label{robin}
\end{equation}
where $\ell(q,p)$ is the Riemannian distance between $p$ and $q$.

The motion of a single vortex depends not only on its initial position but also
on the initial value of a harmonic velocity field (a background flow) \cite{gustafsson2022dipole}.
In the following statement  \cite{ragvil} \cite{errata}
the initial background flow is assumed to be equal to zero:

\nd
{\it A vortex initially placed at any point on a surface $S$ with Riemannian metric $g$
  remains at rest if, and only
  if, the Robin function $R$ associated with  $g$ is constant. A Riemannian metric with
  this property is called 
  a ``Steady Vortex Metric''. }

 The first main result in this paper is the following.
\begin{theorem}[Steady Vortex Metric]
\label{crc}
Let $S$ be a compact Riemann surface.
There exists at least one  steady vortex metric $g$
compatible with the conformal structure of $S$.
There are  examples  where  $g$ is not unique.
\end{theorem}

 The theorem effectively says
  that there always exist a metric for which the Robin function is constant.

This theorem  is a direct consequence of a  theorem proven by K. Okikiolu
\cite{okitorus} \cite{okigenus}   and its proof is
given in Appendix \ref{proofsec}.


The  second theorem in this Section  requires some definitions.
 A one-form $\theta$ on $S$ is  harmonic
if $d\theta=0$ and $d\ast\theta=0$. Since $\ast$ rotates one-forms by
$\pi/2$, harmonic forms are conformal invariants.
The vector space of harmonic forms on $S$ is finite and has dimension $2\calG$
\cite{derham},
 where  $\calG$ is the genus of $S$. Let $\{\theta_1,\ldots\theta_{2\calG}\}$ be an arbitrary
 orthonormal basis of harmonic one-forms in the sense that
 \begin{equation}
 (\theta_j,\theta_k)=\int_S\theta_j\wedge\ast\theta_k=\dt_{jk}.
 \label{ortho}
 \end{equation}
Note: this definition of orthonormality  depends only on the conformal structure.
\begin{theorem}
\label{sg}
Let $(S,g)$ be a compact oriented Riemannian surface  and 
$\sg$ be  the two-form
\[
\sg=\sum_{k=1}^{2\calG} \theta_k\wedge *\theta_{k}\,.
\]

Then the Robin function $R$ is the only solution,
up to an  additive constant,
of the  equation 
\begin{equation}
  \left(\Dt R+\frac{K}{2\pi}-\frac{2}{V}\right)\mu=-\sg\,.
  \label{Req}
\end{equation}
\end{theorem}

 If the genus $\calG$ of $S$ is zero, then $\sg=0$. So a metric on the
sphere is a Steady Vortex Metric if, and only if, it is
of constant curvature.

If the genus of $S$ is greater than zero, then
$\sg$ is the area form of the Bergman metric. The most common
definition of the Bergman metric (\cite{jostgeom} eqs. 1.4.22 and 1.4.23)
uses  a basis of holomorphic
differentials  $\{\omega_1,\ldots,\omega_\calG\}$ 
that satisfy the orthonormality conditions
$\frac{i}{2}\int_S \omega_j\wedge\ov \omega_k=\delta_{jk}$ (here the overbar
denotes complex conjugation).
A  form $\omega_j$  is holomorphic if,  and only if,
$\omega_j=\theta_j+\ast \theta_j$ for some harmonic differential $\theta_j$
(\cite{farkas1992riemann}, Theorem I.3.11). If we define
$\theta_{j+\calG}=\ast \theta_j$,
$j=1,\ldots,\calG$, then the orthogonality condition
for holomorphic differentials implies
the orthogonality condition for harmonic differentials (\ref{ortho}) and
\begin{equation}
  \sigma=\sum_{k=1}^{2\calG} \theta_k\wedge *\theta_{k}=i\,
  \sum_{j=1}^\calG  \omega_j\wedge\ov \omega_j\,,\quad\text{with}
  \quad \int_S\sigma=2 \calG. \label{hol}
\end{equation}
The Bergman metric normalized as $\sigma/(2\calG)$ can also be defined using the
Jacobian variety associated with $S$  (see
\cite{wentworth1991asymptotics},
\cite{jostgeom}, or  equation 1.25 in \cite{fay}). In several references
\cite{wentworth1991asymptotics} \cite{okigenus} \cite{jorge} the normalized 
Bergman metric
$\sigma/(2 \calG)$  is called by the alternative name ``canonical metric''.

Theorem \ref{sg} appeared  in the work of  Okikiolu \cite{okigenus}
(proposition 2.3) as a ``well known'' result related to
the   Arakelov
Green's function (in  Appendix  \ref{proofsec} I give a more self-contained proof
of Theorem \ref{sg} than that in \cite{okigenus}).
The Arakelov Green's function is used in the definition 
of the ``Arakelov metric'' that
is  characterized  by the equation
 (see \cite{jostgeom} eq. 1.4.24):
\begin{equation}
\frac{K_{A}}{2\pi}\mu_A=(2-2\calG)\frac{\sg}{2\calG}\,,\qquad \calG\ge 1\,,\label{arak}
\end{equation}
where: $\mu_A$ and $K_A$ denote the area form and the curvature
of the Arakelov metric.

Equation (\ref{Req}) implies that the several ``natural'' metrics considered in  this paper satisfy
the following relations:
\begin{equation}
  \renewcommand{\arraystretch}{1.7}
  \begin{array}{lll}
  &  \left(\frac{K_{svm}}{2\pi}-\frac{2}{V}\right)\mu_{svm}=-\sg &\quad(\text{SVM})\\
  &  \left(\Dt R_{cc}-\frac{2\calG}{V}\right)\mu_{cc}=-\sg &\quad (\text{constant curvature=CC})\\
  & \left(\Dt R_B+\frac{K_B}{2\pi}\right)\sg=(2-2\calG)\frac{\sg}{2\calG} &\quad
 (\text{Bergman})\\
                  &  \left(\Dt R_A-\frac{2}{V}\right)\mu_A=-\frac{\sg}{\calG}&\quad
    (\text{Arakelov})
    \end{array}\label{rels}
  \end{equation}
  From equations (\ref{arak}) and (\ref{rels}) we obtain
  \begin{equation}\begin{split}
       CC=SVM&\Leftrightarrow Bergman=CC\Leftrightarrow Bergman=SVM\,,\\
       Arakelov=SVM&\Leftrightarrow Bergman=Arakelov\Leftrightarrow
       Arakelov=CC\,.
  \end{split}\label{rel2}
  \end{equation}
For $\calG\ge 1$, therefore,  a constant
curvature metric  is 
 a Steady Vortex Metric if and only if  the Bergman
metric has constant curvature.
For $\calG=1$ this is the case, since the flat metric is the Bergman metric and also
the Arakelov metric.

In any closed  surface $S$ of genus $\calG\ge 2$ the curvature of the Bergman metric
is non positive \cite{jostbook} (theorem 5.5.1).
If the curvature of the Bergman metric $K_B$ is non constant in every
$S$, which as far as I know has not been proved, then
constant curvature metrics will  never be SVM for $\calG\ge 2$.
The last theorem in \cite{lewittes} states that  $K_B(p) = 0$  if and only if $S$
is hyperelliptic and $p$ is one of the $2\calG + 2$ classical
Weierstrass points on $S$. Therefore $K_B$ is not constant in hyperelliptic
 surfaces.
 Since every  surface of genus 2 is hyperelliptic, the 
following
theorem holds.
\begin{theorem}\label{genus2}
  No orientable surface of genus 2 and of constant curvature is a
  Steady Vortex Surface.\end{theorem}

The Gauss-Bonnet theorem implies that the average curvature of the Bergman metric
is  \begin{equation}
  K_{Ba}:= \frac{1}{V_B}\int_S K_{B}\sigma= 2\pi \frac{2-2\calG}{2\calG}\,,\quad
  \text{where}\quad V_B=\int_S \sigma=2\calG\,.\label{KBa}
\end{equation}
We define the deviatoric part $K_{B\delta}$ of $K_B$ as :
\begin{equation}
  K_{B\delta}:= K_B- K_{Ba}\, \quad\text{with}\quad
  \int_S K_{B\delta}\, \sigma=0\,. \label{KB}
\end{equation}
Equation (\ref{rels}) then implies that the Robin function
of the Bergman metric satisfies the simple relation
\begin{equation}
  - \Delta R_B= \frac{1}{2\pi}K_{B\delta}\,.\label{KBd}
\end{equation}
This equation  implies that in any conformal coordinates,
  $\{z, \ov z\}$,
$R_B$ has a simple expression in  terms of the potentials $F_j(z)$  of the
holomorphic differentials
$\omega_j=d F_j=F_j'(z)dz=\partial_zF_j(z)dz$ that appear in the definition
of $\sigma$ in equation (\ref{hol}).
Indeed: $i\,
\sum_{j=1}^\calG  \omega_j\wedge\ov \omega_j=i\,\sum_{j=1}^\calG F_j'(z)\ov{F_j'(z)}
dz\wedge d\ov z=\lambda^2_B\frac{i}{2} dz\wedge d\ov z$, with
$\lambda_B^2=2\sum_{j=1}^\calG F_j'(z)\ov{F_j'(z)}$, 
$\ \Delta= \frac{4}{\lambda^2_B}\partial_z\partial_{\,\ov z}$, 
$\ \lambda^2_B=2\partial_z\partial_{\, \ov z}\sum_{j=1}^\calG F_j(z)\ov{F_j(z)}$, and 
$\ K_B=-\frac{4}{\lambda^2_B}\partial_z\partial_{\,\ov z}\log \lambda_B$; imply
\begin{equation}
  R_B(z, \ov z)=\frac{1}{4\pi}\log\left[\sum_{j=1}^\calG F_j'(z)\ov{F_j'(z)}\right]+ 
 \frac{1-\calG}{2\calG} \sum_{j=1}^\calG F_j(z)\ov{F_j(z)}+constant\,.
 \label{RBcoord}
 \end{equation}

The Riemann sphere 
admits a six-dimensional  group of conformal transformations
(the Moebius group) and a three-dimensional group of isometries.
The Pull-back metric $g_1$ of the round metric $g_0$
by a Moebius transformation that 
 is not an isometry satisfies $g_1=\lb^2g_0$ with $\lb^2\ne 1$ 
 almost everywhere. The Robin function associated to $g_1$ is constant because,
 although different from $g_0$, $g_1$ is
 isometric to $g_0$. This type of ``nonuniqueness'' of a steady vortex
 metric within a conformal class will happen whenever the
 group of diffeomorphisms that preserves the conformal structure is larger than
 the group of isometries.
  Since all spheres with constant curvature are isometric to the round sphere, we
 conclude that $g_0$ is the only
 steady vortex metric modulo isometries. 
The question about the uniqueness of steady vortex metrics on tori will
be postponed to Section \ref{torsec}.

 \section{Generalization to higher dimensions.}

 \label{hdim}

 The definition  of hydrodynamic point vortex is restricted to two dimensions.
There is an analogy between
 vortex and  electric charges in two dimensions \cite{ragvil}. Since
 the theory of electrostatics can be  generalized to higher dimensions,
 electrostatics could be the physical guide to the definition
 of an ``eletrostactic force-free metric'' in dimension $n$.
 The idea although interesting leads to some difficulties, which
 will be discussed in the
 next paragraph,  and it will not be pursued  any further.

  The Green's function $G(q,p)$, solution to
 equation (\ref{compacteq}), can be
 understood as the electric potential due to  a positive
  point charge at
 the point $p$ plus  a  uniform distribution
 of negative  charges.
The Robin function $R(p)$ defined in equation (\ref{robin}) is the
 overall potential energy $G(q,p)$ minus the ``singular potential of the point
 charge'',  $-(2\pi)^{-1}\log \ell(q,p)$, the difference being evaluated
 at $p$. 
 The force upon the point
 charge  is $dR(p)$.
  The most natural definition of Robin function in dimension $n\ge 3$
 would be
 \begin{equation}
\lim_{\ell(q,p)\to 0} \left[
G(q,p)-a_n \ell^{n-2}(q,p)\right],
\label{Rnat}
\end{equation}
where $a_n$ is some constant that depends on $n$.
Unfortunately
the Robin function defined in this way is not a smooth function
 unless additional hypotheses are imposed on the Riemannian metric 
 (see \cite{josthab} for a discussion about this definition
 in the context of the conformal Laplacian).
 Another    way to define the Robin function
 would be  first to compute the force upon a small Riemannian ball of radius
 $\ep$ at $p$ and then to take the limit as $\ep\to 0$ to obtain $dR(p)$.
 This procedure may lead to  quite complicated
 computations as $n$ increases.

 From a mathematical point of view regularity is  the key property of the
 Robin function, which in two dimensions is used in the definitions
 of  vortex motion and  force
 upon an electric charge. In order to define the Robin function 
 in dimension greater than two we will
 regularize the $\dt$-distribution, to do  the computations
 in the regularized setting, and then to take  the limit back to 
 recover the $\dt$-distribution.   In order to do all these limits independently of
 coordinates we use the heat equation. This procedure
 naturally associates the Robin
 function with diffusion and Brownian motion. This association  will be further
 addressed in Section \ref{netsec}.

 Let $(M,g)$ be a compact Riemannian manifold. The heat kernel
 $K:M\times M\times \R_+\to M$ is the fundamental
 solution to the heat equation
 \begin{equation}\label{heateq}
 \left(\frac{\partial}{\partial t}-\Delta_q\right)K(q,p,t)=0, \quad
 \text{with} \quad K(q,p,0)=\dt_p(q).
 \end{equation}
 The initial  condition is understood as a distribution, namely,
 for any $\phi\in\cinfty(M)$
 \[
 \int_MK(q,p,t)\phi(q)\mu(q)\to \phi(p)\quad\text{as}\quad t\to 0_+.
 \]
 The heat kernel is a $\cinfty$ symmetric,  $K(q,p,t)=K(p,q,t)$,
  function.
 Let $0<\lb_1\le\lb_2\le\lb_3\le\ldots$ be
 the nontrivial  eigenvalues  to the problem $\Dt\phi+\lb\phi=0$
 and $\phi_1,\phi_2,\ldots$  be a corresponding 
 $L_2-$orthonormal basis of  eigenfunctions for functions that integrate
 to zero over $M$. Then the 
 the spectral decomposition of the heat kernel is
 \[
 K(q,p,t)=\frac{1}{V}+\sum_{k=1}^\infty\erm^{-\lb_k t}\phi_k(q)\phi_k(p),
 \]
 with pointwise convergence
 (see, for instance,  \cite{sros} for
 basic properties of the heat kernel).
 The Green's function $G(q,p)$ is related to the heat kernel in the
 following way
 \[
 G(q,p)=\int_0^\infty \left( K(q,p,t)-\frac{1}{V}\right)dt
 \]
 This is the formula that allows for the definition of the Robin
 function in dimension $n>2$ by means of the regularization of
 the heat kernel.

 As before, let $\ell(q,p)$ denote the Riemannian distance between
 $q$ and $p$. There exists $\ep>0$ and a set of functions
 $u_0(q,p),u_1(q,p),\ldots$  such that for any given integer $N\ge 0$
 the following estimate holds
 ( the so-called Minakshisundaram-Plejel asymptotic expansion \cite{MP}; see Equations (7)–(9) and the accompanying text)
 \begin{equation}
   \label{asymp1}
 \bigl| K(q,p,t)-\frac{\erm^{\frac{-\ell^2(q,p)}{4t}}}{(4\pi t)^{n/2}}
   \sum_{k=0}^N u_k(q,p)t^k\bigr|
   \le C_{N}t^{N+1-n/2},
   \end{equation}
   for all $(q,p)$ with $\ell(q,p)<\ep$ and  all $t\in(0,1)$,
   where $C_N$ is a constant that  depends only on $N$
   ( see, for instance, \cite{sros} exercise 5 in Section 3.3).
   The functions $u_k$ are  $\cinfty$ and symmetric
   $u_k(q,p)=u_k(p,q)$ \cite{moretti}. If $q=p$ then the
   above expression implies
 \begin{equation}
   \label{asymp2}
 K(p,p,t)= \frac{1}{(4\pi t)^{n/2}}
   \left[a_0(p)+ta_1(p)+\ldots+t^Na_N(p)\right] +E_N(p,t)
   \end{equation}
 where $|E_N(p,t)|< C_Nt^{N+1-n/2}$ for all $p\in M$ and $t\in (0,1)$.
 The functions $a_k(p)=u_k(p,p)$ are local  heat invariants of $M$ that
 can be expressed in terms
   of powers of the Laplacian and the distance function $\ell$
   (\cite{polt}, Theorem 1.2.1). For instance, $a_0(p)=1$ and    $a_1(p)=s(p)/6$,
   where $s(p)$ is the scalar curvature (\cite{sros}, proof of Lemma 3.26 and Proposition 3.29, respectively).

   Suppose that $n\ge 2$ is even and $N$ in equation
   (\ref{asymp2}) is chosen as $\frac{n}{2} -1$.
   Then for $0<\ep<1$ equation (\ref{asymp2}) implies 
  \[
  \int_\ep^1K(p,p,t)dt=
\frac{1}{(4\pi)^{n/2}}
\left[\sum_{k=0}^{\frac{n}{2}-2}a_k(p)
  \frac{\ep^{k+1-\frac{n}{2}}}{\frac{n}{2}-k-1}
  \quad -a_{\frac{n}{2}-1}(p)
  \log\ep
\right] +
C(p,\ep)
\]
where $\lim_{\ep\to 0_+} C(\cdot,\ep)$ is a $\cinfty$ function on $M$.  
 Similarly, if $n\ge 3$ is odd
and $N$ in equation (\ref{asymp2}) is chosen as
$\frac{n}{2}-\frac{3}{2}$ then 
   \[
  \int_\ep^1K(p,p,t)dt=
\frac{1}{(4\pi)^{n/2}}
\left[\sum_{k=0}^{\frac{n}{2}-\frac{3}{2}}a_k(p)\frac{\ep^{k+1-\frac{n}{2}}}
  {\frac{n}{2}-k-1}
\right] +
\tilde C(p,\ep)
\]
where $\lim_{\ep\to 0_+} \tilde C(\cdot,\ep)$ is a $\cinfty$
function on $M$.  
These computations motivate the following definition of the Robin
function $R(p)$. If $n\ge 2$ is even then
\begin{equation}
  \begin{split}
    R(p)=\lim_{\ep\to 0_+}\Biggl\{&
    \int_\ep^\infty\left(K(p,p,t)-\frac{1}{V}\right)dt  \\
    &-
    \frac{1}{(4\pi)^{n/2}}
    \left[\sum_{k=0}^{\frac{n}{2}-2}a_k(p)
      \frac{\ep^{k+1-\frac{n}{2}}}{\frac{n}{2}-k-1} \
      -a_{\frac{n}{2}-1}(p)\big(\log(4\ep)-\gamma\big)\right]
\Biggr\}\,,
\end{split}
\label{Re}
\end{equation}
where
$\gamma=-\int_0^\infty \erm^{-x}\log x\, dx=0.577215\ldots$ is the Euler's constant.
The constant term $a_{\frac{n}{2}-1}(p)(\log 4-\gamma)/(4\pi)^{n/2}$ were added to the
right-hand side of equation (\ref{Re}) to preserve the definition
of the Robin function given in equation (\ref{robin}).

   If $n\ge 3$ is odd then
   \begin{equation}
     \begin{split}
  R(p)=\lim_{\ep\to 0_+}\Biggl\{&
  \int_\ep^\infty\left(K(p,p,t)-\frac{1}{V}\right)dt\\
  &-
\frac{1}{(4\pi)^{n/2}}
\left[\sum_{k=0}^{\frac{n}{2}-\frac{3}{2}}a_k(p)
  \frac{\ep^{k+1-\frac{n}{2}}}{\frac{n}{2}-k-1}\right]
\Biggr\}.
\end{split}
\label{Ro}
\end{equation}

   \begin{theorem}
     \label{Rth}
     The Robin function can be written in the
     following alternative way: for $n\ge 3$ odd
\[
       R(p)=\lim_{q\to p}\left\{
        G(q,p)-\frac{1}{(4\pi)^{n/2}}
   \sum_{k=0}^{\frac{n}{2}-\frac{3}{2}}u_k(q,p)
\left(\frac{4}{\ell^2(q,p)}\right)^{\frac{n}{2}-k-1}
      \Gamma(\frac{n}{2}-k-1)\right\}
\]
and for $n\ge 2$ even
\[
\begin{split}
  R(p)=&\lim_{q\to p}\Biggl\{ 
   G(q,p)-
    \frac{1}{(4\pi)^{n/2}}
    \sum_{k=0}^{\frac{n}{2}-2}u_k(q,p)
  \left(\frac{4}{\ell^2(q,p)}\right)^{\frac{n}{2}-k-1}\Gamma(\frac{n}{2}-k-1)
      \\
      &\qquad 
      +
    \frac{1}{(4\pi)^{n/2}}
      u_{\frac{n}{2}-1}(q,p)\Bigl[\log\ell^2(q,p)
      \Bigr]
    \Biggr\}\,.
    \end{split}
\]
\end{theorem}   
  \bproof At first we prove the statement for  $n$ odd.
   From equation (\ref{Ro})
   \begin{eqnarray}
  R(p)&=&\lim_{\ep\to 0_+}\lim_{q\to p}\Biggl\{
  \int_0^\infty\left(K(q,p,t)-\frac{1}{V}\right)dt
  -\int_0^\ep\left(K(q,p,t)-\frac{1}{V}\right)dt\nonumber\\
  &&-
\frac{1}{(4\pi)^{n/2}}\left.
\left[\sum_{k=0}^{\frac{n}{2}-\frac{3}{2}}u_k(q,p)
  \frac{\ep^{k+1-\frac{n}{2}}}{\frac{n}{2}-k-1}\right]
\right\}\label{lqpa}
\end{eqnarray}
   The term $\int_0^\ep K(q,p,t)dt$ is estimated in the following way.
   From  equation 
   (\ref{asymp1})
   \[
   K(q,p,t)=\frac{\erm^{\frac{-\ell^2(q,p)}{4t}}}{(4\pi t)^{n/2}}
   \sum_{k=0}^{\frac{n}{2}-\frac{3}{2}} u_k(q,p)t^k+ \Phi_1(q,p,t)
   \]
   such that $|\Phi_1(p,p,t)|\le C_1 t^{-1/2}$, $C_1>0$,
   for $t\in(0,1)$ (the constants $C_1,C_2\ldots$ do not depend on
   $t$, $q$, $p$, or $\ep$).
   Therefore
   \[
   \int_0^\ep K(q,p,t)dt=\frac{1}{(4\pi)^{n/2}}
   \sum_{k=0}^{\frac{n}{2}-\frac{3}{2}}u_k(q,p)
   \int_0^\ep \erm^{\frac{-\ell^2(q,p)}{4t}}
    t^{-\frac{n}{2}+k}dt+ \Phi_2(q,p,\ep)
    \]
    such that $|\Phi_2(p,p,\ep)|<C_2\ep^{1/2}$, $C_2>0$.
    With the change of variables  $t=\frac{\ell^2}{4s}$ we obtain
    \[
    \int_0^\ep \erm^{\frac{-\ell^2}{4t}}
    t^{-\frac{n}{2}+k}dt=\left(\frac{\ell^2}{4}\right)^{k+1-\frac{n}{2}}\left[
    \int_0^\infty\erm^{-s}s^{\frac{n}{2}-k-2}ds-
    \int_0^{\frac{\ell^2}{4\ep}}\erm^{-s}s^{\frac{n}{2}-k-2}ds\right]
    \]
    where $\int_0^\infty\erm^{-s}s^{\frac{n}{2}-k-2}ds=\Gamma(\frac{n}{2}-k-1)$
    is the Gamma function. Using that $\erm^{-s}=1-sF(s)$, where
    $F(s)=\int_0^1\erm^{-\eta s} d\eta$, an explicit computation gives
      \[
      \left(\frac{\ell^2}{4}\right)^{k+1-\frac{n}{2}}
      \int_0^{\frac{\ell^2}{4\ep}}\erm^{-s}s^{\frac{n}{2}-k-2}ds=
      \frac{\ep^{-\frac{n}{2}+k+1}}{\frac{n}{2}-k-1}-\frac{\ell^2}{4}
      \Phi_3\left(\frac{\ell^2}{4\ep}\right)
      \]
      where $|\Phi_3\left(\frac{\ell^2}{4\ep}\right)|<1$.
      Therefore,
      \[
      \int_0^\ep \erm^{\frac{-\ell^2}{4t}}
      t^{-\frac{n}{2}+k}dt=\left(\frac{\ell^2}{4}\right)^{k+1-\frac{n}{2}}
      \Gamma(\frac{n}{2}-k-1)-
      \frac{\ep^{-\frac{n}{2}+k+1}}{\frac{n}{2}-k-1}+\frac{\ell^2}{4}
      \Phi_3\left(\frac{\ell^2}{4\ep}\right)
      \]
      Finally, using that
      $\int_0^\infty\left(K(q,p,t)-\frac{1}{V}\right)dt=G(q,p)$ and
      substituting all the previous estimates into equation
      (\ref{lqpa}) we obtain
      \begin{eqnarray*}
        R(p)&=&\lim_{\ep\to 0_+}\lim_{q\to p}\left\{
        G(q,p)-\frac{1}{(4\pi)^{n/2}}
   \sum_{k=0}^{\frac{n}{2}-\frac{3}{2}}u_k(q,p)
\left(\frac{\ell^2}{4}\right)^{k+1-\frac{n}{2}}
      \Gamma(\frac{n}{2}-k-1)\right.\\
      &&\left. +\frac{\ep}{V}-\Phi_2(q,p,\ep)-\frac{\ell^2(q,p)}{4}
      \frac{1}{(4\pi)^{n/2}}
      \sum_{k=0}^{\frac{n}{2}-\frac{3}{2}}u_k(q,p)
      \Phi_3\left(\frac{\ell^2(q,p)}{4\ep}\right)\right\}.
      \end{eqnarray*}
      The limit  as $q\to p$ of the second line of this equation is
      $\frac{\ep}{V}$. Since for a fixed $\ep$, the limit
      as $q\to p$ of the
      expression inside brackets exists then the limit as $q\to p$ 
      of the sum in the first line also exists and  does not depend
      on $\ep$. So, the proof for $n$ odd is finished.

      The proof for  $n$ even is similar. The only difference is
      that it is necessary to estimate the additional integral
    \[
    \int_0^\ep \erm^{\frac{-\ell^2}{4t}}
    t^{-1}dt=\int_{\frac{\ell^2}{4\ep}}^\infty\erm^{-s}s^{-1}ds=
    -\log\left(\frac{\ell^2}{4\ep}\right)+
    \int_{\frac{\ell^2}{4\ep}}^\infty\erm^{-s}\log s\,  ds.
    \]
    The last integral  is equal to minus the Euler's constant  as $q\to p$. 
   \eproof

    We remark that for $n>2$
    the term  of highest order in $\ell^{-2}$ is
\[
\frac{1}{(4\pi)^{n/2}}
\left(\frac{4^{\frac{n}{2}-1}}{\ell^{n-2}(q,p)}\right)\Gamma(\frac{n}{2}-1)\,,
\]
where we used that $u_0(p,p)=a_0(p)=1$, is minus
the ``Newtonian potential'' that appears
in equation (\ref{Rnat}).

For $n=2$, theorem \ref{Rth} states that the Robin function defined by
equation (\ref{Re}) coincides with that
given in equation (\ref{robin}). 

Theorem \ref{Rth} can also  be obtained from the Hadamard parametrix,
see \cite{garabedian1986} section 5.3.

The Robin function as given in Theorem \ref{Rth}
can be written in terms of the analytic extension of the 
Minakshisundaram–Pleijel zeta function, 
\cite{steiner2005geometrical} (dimension two)
and \cite{bilal2013multi} (dimension greater than one). The relation between the
Robin function and the zeta-function is presented in Appendix \ref{min}.

   \section{The ``Narrow Escape Time (NET)''.}
\label{netsec}

In the context of a compact boundaryless manifold $M$, the narrow escape
problem can be described in the following way.
Consider a Brownian motion  on $M$, whose infinitesimal
generator is the  Laplace-Beltrami operator $\Delta$.
Let
$B_\epsilon(q)\subset M $ be a geodesic ball of small radius $\epsilon>0$.
This ball will be the absorbing set or  the
small window through which a particle  can escape.
The  amount of time that a  particle initially at $p$ is expected to spend in
$ M\backslash B_\epsilon(q)$ (the mean sojourn time) will be denoted as
$v_\epsilon(p,q)$. This function  is the ``narrow escape time'' (NET)
since
 it measures the mean time it takes for a particle initially at $p$ to escape
 through the narrow window $B_\epsilon(q)$.
  The NET is the solution to the problem
  (see  \cite{holcman2014narrow}, equation 3.1):
\begin{equation}\label{Greenep}\begin{split}
    & \nu  \Delta_p v_\epsilon(p,q)=-1\,, \quad  p\in M\backslash B_\epsilon(q)\,,\\
    & \text{with}\quad
     v_\epsilon(p,q)=0\quad\text{for}
     \quad p\in \partial B_\epsilon(q)\,,
     \end{split}
   \end{equation}
   where $\nu$ is a diffusion coefficient with dimensional units length$^2/$time.
   The NET averaged against a uniform distribution of initial points in
   $M\backslash B_\epsilon(q)$,
   \begin{equation}
     \ov v_\epsilon(q)=\frac{1}{V-|B_\epsilon(q)|}
     \int_{M\backslash B_\epsilon(q)}
     v_\epsilon(p,q)\mu(p),\label{ovv}
   \end{equation}
   gives the expected time a particle randomly placed in the manifold remains in it
   until it scapes through $ B_\epsilon(q)$.

    In dimension 2 the following theorem was proved in \cite{doyle2017spectral}
   (Lemma 4.1 and Theorem 4.2 part 2).
   \begin{theorem}
     \label{tnet} In dimensions 2, 3, and 4,
     the ``Narrow Escape Time'' (NET) is given by
     \begin{equation}
       v_\epsilon(p,q)=-\frac{V}{\nu}\, G(p,q) +  \ov v_\epsilon (q)+E(p,q,\epsilon)\,,
       \label{net1}
     \end{equation}
where $\lim_{\epsilon\to 0}E(p,q,\epsilon)=0$.
     The average NET,  equation {\rm (\ref{ovv})}, is given by
     \[ \renewcommand{\arraystretch}{1.7}
       \begin{array}{rll}
         \ov v_\epsilon (q)&=\frac{V}{\nu}
            \left\{-\frac{1}{2\pi}\log\epsilon+R(q)+E_2(q,\epsilon)\right\}
                              &n=2\\
    \ov v_\epsilon (q)&=\frac{V}{\nu}
            \left\{\frac{1}{4\pi} \epsilon^{-1}+R(q)+E_3(q,\epsilon)\right\}
                              & n=3\\
          \ov v_\epsilon (q)&=\frac{V}{\nu}
            \left\{\frac{1}{4\pi^2} \epsilon^{-2}
   -\frac{1}{48\pi^2}S(q)\log\epsilon+\frac{1}{192\pi^2}S(q)
                              +R(q)+E_4(q,\epsilon)\right\}
                              &n=4\,,
        \end{array}
      \]
      where $S(q)$ is the scalar curvature at $q$ and
      $\lim_{\epsilon\to 0}E_n(p,q,\epsilon)=0$.
      \end{theorem}

      \nd {\it Remarks:}

      \nd$\bullet$ The normalization $\int_M  G(p,q)\mu(p)=0$ (equation (\ref{G}))
      ensures the compatibility of both sides of equation (\ref{net1}).

      \nd$\bullet$ The  NET increases as $\epsilon$  decreases in the same way as
      the Newtonian potential in  $\R^n$ 
      increases as the distance to the singularity decreases (see e.g.
      \cite{holcman2014narrow}, section 3, for the same result for surfaces).
      This is true in all dimensions, not only $n=2,3,4$.
      
      \nd $\bullet$
      In dimensions 2 and 3,   the divergent terms of  $\ov v_\epsilon(q)$
     with respect to  $\epsilon$ do not depend on $q$.
      For $n=4$ this is no longer  true, since $\ov v_\epsilon (q)$
      contains a logarithmic
      divergent term that is proportional to the mean curvature $S(q)$. 
      If the mean curvature is constant on $M$, then  the
      dependence of $\ov v_\epsilon(q)$ on $q$ as  $\epsilon\to 0$
      is   determined by the Robin function, 
      as it is in dimensions
      2 and 3.


      \vskip 2truemm
      \bproof We will prove only the case n=4. The proof of the cases n=2 and
      n=3 is simpler and goes along  the same lines.

      We write $v_\epsilon(p,q)=-\frac{V}{\nu}\, G(p,q) + \frac{V}{\nu}
h_\epsilon(p,q)$ and from equations
(\ref{compacteq})
  and   (\ref{Greenep}) we obtain
      \begin{equation}
     \Delta_p h_\epsilon(p,q)=0\,, \quad  p\in M\backslash B_\epsilon(q)\,, \quad \text{with}\quad
     h_\epsilon(p,q)= G(p,q)\quad\text{for}
     \quad p\in \partial B_\epsilon(q)\,.
     \label{heps1}
     \end{equation}

     Let $x\in \R^n$ be an orthonormal coordinate system on
     the tangent space of $M$
     at $q$. Let $x(p):=\exp^{-1}_qp$ be  geodesic normal coordinates in $M$ defined
     in a neighborhood of $q$. The metric tensor in this
     coordinates is given by
     $g_{ij}(x)=\delta_{ij}-\frac{1}{3}R_{ikjl}x^kx^l+\Oc(|x|^3)$.    
     Theorem \ref{Rth} implies that for $p$ sufficiently close to $q$
        \[G(p,q)=
    \frac{u_0(p,q)}{(4\pi)^{2}}
      \left(\frac{4}{|x|^2}\right)
          -
    \frac{ u_{1}(p,q)}{(4\pi)^{2}}
      \log|x|^2
      + R(q)+ {\cal R}_1(p,q)\,.\]
From $u_0(p,q)=1/\sqrt{{\rm det}(\exp_qp)}$ (\cite{sros} equation (3.11))
    we obtain
    \[\begin{split}
      u_0(q,p)&=\frac{1}{({\rm det}\, g)^{1/4}}=1-\frac{1}{12}R_{kil}^i(q)x^k x^l
   +\Oc(|x|^3)\\ &=
   1-\frac{1}{48}S(q)|x|^2-\frac{1}{12}Z_{kl}(q)x^k x^l +
   \Oc(|x|^3) \end{split} \]
    where:
    $R_{jkl}^i$ is the Riemman curvature tensor and $Z_{kl}$ is the traceless
    Ricci tensor. From  \cite{sros} Proposition 3.29,
    \[
      u_1(p,q)=u_1(q,q)+\Oc(|x|)=\frac{1}{6}S(q)+\Oc(|x|)\,.\]
    
    The expressions in the previous paragraph imply
    that for $p$ sufficiently close to $q$  
   \begin{equation}
     G(p,q)=\underbrace{\frac{1}{4\pi^2|x|^2}-\frac{S(q)}{192 \pi^2}
     -\frac{S(q)}{48 \pi^2}\log|x|
           + R(q)}_{\text{term a}}-
     \underbrace{\frac{ Z_{kl}(q)}{48 \pi^2}
     \frac{x^k x^l}{|x|^2}}_{\text{term b}}+ \underbrace{{\cal R}_2(p,q)}_{\text{term c}}\,,\label{Gest}
     \end{equation}
         where  $\lim_{|x|\to 0}{\cal R}_2(p,q)=0$.

         The solution $h_\epsilon$ to the problem in equation (\ref{heps1})
         can be split into three terms a, b, and c, according to the
         decomposition of the boundary conditions as given in equation
         (\ref{Gest}).
 The  term a is constant
         for $|x|=\epsilon$. This  
         term  appears in the expression for $\ov v_\epsilon (q)$ in the statement
         of the theorem.

         The maximum principle dictates that the maximum of the function
         $p \mapsto |h_c(p,q,\epsilon)|$, where $h_c(p,q,\epsilon)$ is the solution to 
  \[
     \Delta_p h_{c,\epsilon}(p,q) = 0, \quad p \in M \setminus B_\epsilon(q), \quad \text{with} \quad
     h_{c,\epsilon}(p,q) = {\cal R}_2(p,q) \quad \text{for} \quad
     p \in \partial B_\epsilon(q),
  \]
   is attained on $\partial B_\epsilon(q)$. Given that
   $h_{c,\epsilon}(p,q) = {\cal R}_2(p,q)$ for $p \in \partial B_\epsilon(q)$, the limit $\epsilon \to 0$ implies $p \to q$, and 
   since $\lim_{p \to q}{\cal R}_2(p,q) = 0$, it follows that $\lim_{\epsilon \to 0}|h_c(p,q,\epsilon)| \to 0$
   for $p \in M \setminus B_\epsilon(q)$. Consequently, the component of $h_\epsilon$ corresponding to term (c)
   in equation (\ref{Gest}) contributes to the function $E(q, p, \epsilon)$ as stated in the theorem.

         The part of $h_\epsilon$ associated  the term b in equation (\ref{Gest}),
           will be   denoted as $H_\epsilon(p,q)$. It satisfies
         the problem 
 \begin{equation}
     \Delta_p H_\epsilon(p,q)=0 \quad \text{with}\quad
     H_\epsilon(p,q)= -
     \frac{ Z_{kl}(q)}{192 \pi^2}
    \frac{x^k x^l}{\epsilon^2}\quad\text{for}
     \quad |x|=\epsilon\,.
     \label{h2}
     \end{equation}
In order to finish the proof we must show that $|H_\epsilon(p,q)|\to 0$
     as $\epsilon\to 0$ $|H_\epsilon(p,q)|\to 0$
     as $\epsilon\to 0$Y. The proof has several steps.
\begin{proposition}\label{av}
       \[
         \int_{|x|=\epsilon} Z_{kl}(q)
           \frac{x^k x^l}{\epsilon^2}d\sigma(x)=\Oc(\epsilon^{5})\,,
         \]
         where $d\sigma(x)$ is the ``area form'' on the geodesic
         sphere $\partial B_\epsilon(q)$.
       \end{proposition}
       \bproof The area form on  $\partial B_\epsilon(q)$ satisfies
       $d\sigma(x)=[1-\frac{1}{6}R^i_{kil}x^kx^l+\Oc(|x|^3)]d\sigma_{E}(x)$,
       where $d\sigma_{E}(x)$ is the Euclidean area form on the sphere
       $|x|=\epsilon$. The function $ Z_{kl}(q)x^k x^l$ is harmonic with respect
       to the Euclidean Laplacian, since the trace of $Z$ is zero, and therefore
       its integral over  $|x|=\epsilon$ with respect to  $d\sigma_{E}(x)$
        is zero.
       The proposition follows from the expression for $d\sigma(x)$ and
       $|d\sigma_E|=\Oc(\epsilon^{3})$ on $|x|=\epsilon$.
       \eproof

       The identity $\Delta_z[G(z,q)-G(z,p)]=\delta_p(z)$,  for $z$ and $p$ in
       $M\backslash B_{2\epsilon} (q)$, and Green's second identity
       imply that for $\epsilon$ sufficiently small
       \begin{equation}\begin{split}
         H_\epsilon(p,q)=&\int_{\partial B_{2\epsilon}}H_\epsilon(x,q)
         \nabla_x[G(x,q)-G(x,p)]\cdot\frac{x}{|x|}d\sigma(x)\\
         &-
         \int_{\partial B_{2\epsilon}}[G(x,q)-G(x,p)]
         \nabla_x H_\epsilon(x,q)\cdot\frac{x}{|x|}d\sigma(x)\,.
         \end{split}\label{eqHe}
       \end{equation}

 We will first estimate  the integral in the second line of equation (\ref{eqHe}).
 Equation (\ref{Gest}) implies that $G(x,q)-G(x,p)$ with $|x|=2\epsilon$ can
 be written
 as a term $A_1=\Oc(\epsilon^2)$ that does not depend on $x$ and
 a term $A_2$ that
 is bounded by a constant $C_1(q)$ that is
 independent of $x$ and $\epsilon$.
 The integral $ A_1\int_{\partial B_{2\epsilon}}
 \nabla_x H_\epsilon(x,q)\cdot\frac{x}{|x|}d\sigma(x)=0$
 because 
$\int_{M\backslash B_{2 \epsilon}(q)}\Delta_p H_\epsilon(p,q)\mu(p)=0$.
In order to estimate the integral that contains $A_2$ we will use
one of the Schauder interior estimates
\cite{gilbarg2001elliptic} (Corollary 6.3)\footnote{Here is the reason for
  having integrated over the domain $M\backslash B_{2\epsilon} (q)$ and not
     $M\backslash B_{\epsilon} (q)$.
}    
\[
 \epsilon \max_{|x|=2\epsilon}\left| \nabla_x H_\epsilon(z,q)\cdot\frac{x}{|x|}\right|\le
  \max_{p\in M\backslash B_{\epsilon} (q)}\left| H_\epsilon(p,q)\right|
  \le
  \max_{|x|=\epsilon}\left|H_\epsilon(x,q)\right|=C_2(q)\,,\]
where the second inequality follows from  the maximum principle and the constant
$C_2(q)$ does not depend on $\epsilon$. So
$\left|\int_{|x|=2\epsilon} A_2
  \nabla_x H_\epsilon(x,q)\cdot\frac{x}{|x|}d\sigma(x)\right|= \Oc(\epsilon^2)$ and
the integral in the second line of equation (\ref{eqHe}) is at most of the 
order of $\epsilon^2$. It remains to
estimate   the integral in the first line of equation (\ref{eqHe}).

For a fixed $p\in M\backslash B_{2\epsilon} (q)$ the function 
$\nabla_x[G(x,p)]\cdot\frac{x}{|x|}$ restricted to  $|x|=2\epsilon$ is
uniformly bounded with respect to $\epsilon$.  Therefore, using that
 $|H_\epsilon(x,q)|<C_2(q)$, we obtain 
$\int_{|x|=2\epsilon}H_\epsilon(x,q)
\nabla_x[G(x,p)]\cdot\frac{x}{|x|}d\sigma(x)=\Oc(\epsilon^3)$ and
it remains to estimate
$\int_{\partial B_{2\epsilon}}H_\epsilon(x,q)
\nabla_x[G(x,q)]\cdot\frac{x}{|x|}d\sigma(x)$.

The  term $\frac{1}{4\pi^2|x|^2}$
in  equation (\ref{Gest}) is the leading order term of a  parametrix for
the Laplace equation (see \cite{garabedian}, equation (5.79), or
\cite{aubin}, theorem 4.13 equation (17)). This implies that
$G(x,q)-\frac{1}{4\pi^2|x|^2}$, where $G(x,q)$ is given in equation
(\ref{Gest}),  can be differentiated for  $x\ne0$ and the derivative of
${\cal R}_2(x,q)$
is dominated by  those of the other terms,
so that 
$|\nabla_x [G(x,q)-\frac{1}{4\pi^2|x|^2}]\cdot\frac{x}{|x|}|=\Oc(1/|x|)$.
This and   $|H_\epsilon(x,q)|<C_2(q)$ 
imply
\[\int_{|x|=2\epsilon}H_\epsilon(x,q)
\nabla_xG(x,q)\cdot\frac{x}{|x|}d\sigma=
\int_{|x|=2\epsilon}H_\epsilon(x,q)
\nabla_x\left[\frac{1}{4\pi^2|x|^2}\right]
\cdot\frac{x}{|x|}d\sigma+\Oc(\epsilon^2)\,.\]
It remains  to estimate the integral in the right-hand side of
this equation.

  Green's second identity with  $\Delta_x H_\epsilon(x,q)=0$
and 
$\nabla_x\left[\frac{1}{4\pi^2|x|^2}\right]\cdot\frac{x}{|x|}=-\frac{1}{2\pi^2|x|^3}$, which is valid because $x$ are normal coordinates,
imply
\[\begin{split}
 & \int_{B_{2\epsilon}(q)\backslash B_{\epsilon}(q)} H_\epsilon(x,q)\Delta
    \left[\frac{1}{4\pi^2|x|^2}\right]dx^4=\\
    &
    -\frac{1}{16\pi^2\epsilon^3} \int_{|x|=2\epsilon}
    H_\epsilon(x,q) d\sigma
     +\frac{1}{2\pi^2\epsilon^3} \int_{|x|=\epsilon}
    H_\epsilon(x,q) d\sigma\\
&    -
  \frac{1}{16\pi^2|\epsilon|^2}  \int_{|x|=2\epsilon}
  \nabla H_\epsilon(x,q)\cdot \frac{x}{|x|} d\sigma
  +  \frac{1}{4\pi^2|\epsilon|^2}  \int_{|x|=\epsilon}
  \nabla H_\epsilon(x,q)\cdot \frac{x}{|x|} d\sigma\,.
  \end{split}
\]
The integrals in the last line  are zero because 
$\int_{M\backslash B_{s \epsilon}(q)}\Delta_p H_\epsilon(p,q)\mu(p)=0$, for $s=1,2$.
Due to equation (\ref{h2}) and  Proposition \ref{av},
$\frac{1}{2\pi^2\epsilon^3} \int_{|x|=\epsilon}
H_\epsilon(x,q) d\sigma=\Oc(\epsilon^2)$. A computation using the
expression for the Laplacian in geodesic normal coordinates
\cite{sros} (Theorem  2.63) gives $\Delta
\left[\frac{1}{4\pi^2|x|^2}\right]=\Oc(|x|^{-2})$.
This and  $|H_\epsilon(x,q)|<C_2(q)$ imply
$\int_{B_{2\epsilon}(q)\backslash B_{\epsilon}(q)} H_\epsilon(x,q)\Delta
\left[\frac{1}{4\pi^2|x|^2}\right]dx^4=\Oc(\epsilon^2)$.
In conclusion, all these estimates imply
\[\int_{|x|=2\epsilon}H_\epsilon(x,q)
\nabla_x\left[\frac{1}{4\pi^2|x|^2}\right]
\cdot\frac{x}{|x|}d\sigma=-
\frac{1}{16\pi^2\epsilon^3} \int_{|x|=2\epsilon}
H_\epsilon(x,q) d\sigma=\Oc (\epsilon^2)\,,\]
which finishes the proof.
\eproof

   \section{ Examples of  non-constant curvature  uniform drainage
    surfaces: Okikiolu's tori.} 
\label{torsec}

The flat metric $g_0$ on any two-dimensional torus is a steady vortex metric (SVM).
Equation (\ref{Req}) implies that there exists a second
SVM $g_1$ conformal to $g_0$, $g_1=\lb^2g_0$, if and only if 
\begin{equation}
  \left(\frac{K_1}{2\pi}-\frac{2}{V_1}\right)\mu_1=
   -\frac{2}{V_0}\mu_0
  \label{Req2}
\end{equation}
Normalizing the volumes $\mu_0$ and $\mu_1$ such that $V_0=V_1=1$, using
$-\Dt_0\log\lambda=\lambda^2 K_1$ and  $\mu_1=\lambda^2 \mu_0$, and
defining $f=\log\lb^2$ 
we get the  following equation for $f$
\begin{equation}
\Dt_0f=8\pi-8\pi \erm^f \label{hunique}
\end{equation}
To each nontrivial solution to this equation corresponds  a SVM  $g_1$
conformal to $g_0$.

In the following we present a family of examples due to  Okikiolu \cite{okitorus}
of  non flat  2-dimensional tori that have  constant Robin function, and
so are   uniform drainage surfaces. Each non flat torus in the family
is conformal to a flat torus, which is also a uniform drainage surface. The
Robin function of the non flat tori are smaller than those of the
conformally equivalent flat tori, and  so the narrow escape time of the
non flat tori are smaller than those of the conformally equivalent flat tori.
 There are two differences
between our presentation and that of Okikiolu. We simplify the 
proof that the Robin functions of the non flat tori are smaller
than those of the flat tori and   we represent  the non flat tori in
$\R^3$ as the
quotient of an isometrically embedded cylinder.

Consider the torus $\R^2/(a\Z\times a^{-1}\Z)$, $a\ge 1$,
with the conformal structure induced by the flat metric $g_0$.
If $a\le 2/\sqrt \pi$, then  $g_0$ is the unique uniform drainage metric
\cite{nolasco},  and  if $a>\sqrt{\pi/2}$, then $g_0$ is not unique
\cite{lin2006uniqueness}.
When 
$a>\sqrt{\pi/2}$ a second natural vortex metric can be constructed 
in the following way \cite{okitorus}.
Let $(x,y)$ be Cartesian 
coordinates on $\R^2$. We will look for a nontrivial solution to equation 
(\ref{hunique}) that depends only on the variable $x$, $\partial_y f=0$,
 with  $f(x+a)=f(x)$.
Then $f$ must satisfy $\ddot f:=\frac{d^2f}{dx^2}=8\pi(1-\erm^f)$. This 
ordinary differential equation has a single equilibrium and 
a first integral
\begin{equation}
  H(f,p)=p^2/2 +8\pi(\erm^f-f-1)\,,\quad p=\dot{f}\,.\label{Hf}
\end{equation}
This shows that 
  all solutions $f$ are periodic with a period $T(E)$, where $E$ is the value
  of the first integral associated to the solution.
  The linearized period at $(f,\dot f)=(0,0)$ is  $T(0)=\sqrt{\pi/2}$.

  The period  function $E\to T(E)$ of equation
  $\ddot f=8\pi(1-\erm^f)$ was
  studied in
\cite{chicone1987monotonicity} (p. 315), where it is shown that
$\frac{d}{dE}T (E)>0$. We will additionally
show that  $\lim_{E\to\infty} T(E)=\infty$.
   Consider the solution associated to the initial condition $f(0)=0$,
   $\dot f(0)=-\sqrt{2E}$ and integrate the equation $\ddot f=8\pi(1-\erm^f)$
   on the interval
  $[0,\beta]$,  where $\beta>0$ is the smallest value such that
  $f(\beta)=0$. Since $\dot f(\beta)=\sqrt{2E}$,  the result is
  \begin{equation}
    \sqrt{2E}/(4\pi)= \beta-\int_0^\beta\erm^f dt<\beta<T(E)\,, \label{Taux}
  \end{equation}
  and therefore
  $\lim_{E\to\infty} T(E)=\infty$. 
As a result, equation (\ref{hunique}) has  nontrivial solutions for all
$a>\sqrt{\pi/2}$ such that $f(x+a)=f(x)$
(indeed as many different solutions as we wish provided $a$ is
sufficiently large). 

For a given $a>\sqrt{\pi/2}$, let $g_1=\erm^{f(x)}(dx^2+dy^2)$ be the metric
associated to a periodic solution to  $\ddot f=8\pi(1-\erm^f)$ with minimal
period $a$.
We will use lemma \ref{usmooth} to
show that the  Robin function $R_1$ associated to $g_1$
has a smaller value than the Robin function $R_0$ of the flat metric.
The area form associated to  $g_1$ is given by
$\mu_1=\erm^f \mu_0=\left[1-\frac{\ddot f}{8\pi}\right]dx\wedge dy$ and
the equation that
determines the function $\phi$ in  lemma \ref{usmooth}
becomes
\[
  \Delta_0 \phi \ \ dx\wedge dy=\mu_1-\mu_0=-\frac{\ddot f}{8\pi} dx\wedge dy\,,
  \qquad \int_S\phi\mu_0=0\]
that implies
\[
  \phi(x)=-\frac{f(x)}{8\pi} +\frac{1}{8\pi a}\int_0^af(x)dx\,.
\]
The constant $c=
-\frac{1}{V}\int_S\phi(p)\mu_1(p)$ in lemma \ref{usmooth} can be easily
computed and is equal  to
$c=\frac{1}{(8\pi)^2 a}\int_0^a \dot f^2 dx$. These computations and equation
(\ref{t4}) imply
\begin{equation}
  R_1-R_0=\frac{1}{4\pi a}\int_0^afdx
  +\frac{1}{(8\pi)^2 a}\int_0^a \dot f^2  dx
\label{Robindif}
\end{equation}
If we use the first integral $H$ in equation (\ref{Hf}) to eliminate
$f$ in the right-hand side of this equation  and then use
$\frac{1}{a}\int_0^a \erm^fdx=1$, which we obtain integrating 
 $\ddot f=8\pi(1-\erm^f)$ over the interval
 $[0,a]$, then
\begin{equation}
  R_1-R_0=-\frac{H}{32\pi^2}
  +\frac{1}{32\pi^2 a}\int_0^a \dot f^2  dx\,.
  \label{t42}
\end{equation}
 The equation
$\ddot f=8\pi(1-\erm^f)$ can be written in Hamiltonian form with Hamiltonian
function $H$. 
Using the definition
of the action $I(E)=\frac{1}{2\pi}\oint p df$
from   Hamiltonian mechanics \cite{arnol2013mathematical}
we can write
\[
  \frac{1}{32\pi^2 a}\int_0^a \dot f^2  dx=
  \frac{1}{32\pi^2 a}\int_0^a p \dot f  dx=
  \frac{1}{16\pi a}\frac{1}{2\pi}\oint p  df=\frac{I(E)}{16\pi a}
\]
In this expression $a$ is the period of $f$, and therefore $a=T(E)$ where
$E$ is the value of $H$ associated to $f$. The Hamiltonian function
can be written as a function of the action $E=H(I)$ with
$H^\prime(I)=2\pi/T(I)$. All these results imply that
equation (\ref{t42}) can be written
as
\begin{equation}
  R_1-R_0=\frac{1}{32\pi^2}\left(I H^\prime(I)-H(I)\right)
    \label{t43}
\end{equation}
 
\begin{figure}[hptb!]
\centering
\begin{minipage}{0.5\textwidth}
\centering
 \includegraphics[width=0.9\textwidth]{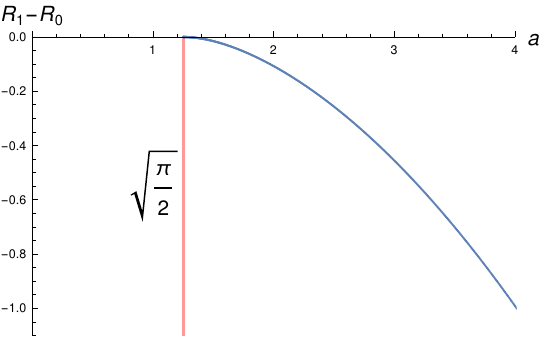}
\end{minipage}\hfill
\begin{minipage}{0.5\textwidth}
\centering
\includegraphics[width=0.9\textwidth]{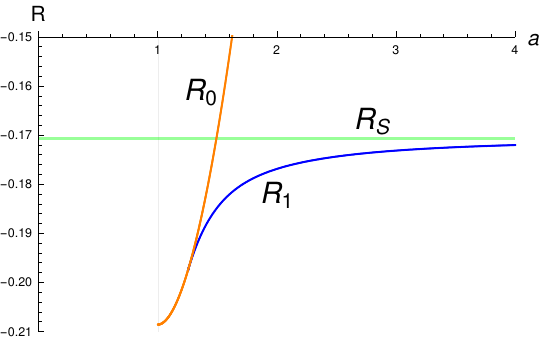}
\end{minipage}\hfill
\caption{LEFT:
  Difference $R_1-R_0$ as a function of $a$, where
  $R_1$ ($R_0$) is  the Robin function
  of the non flat torus $\{ \R^2/(a\Z\times a^{-1}\Z),g_1\}$
  (flat torus  $\{ \R^2/(a\Z\times a^{-1}\Z),g_0\}$).
 RIGHT: Graphs of $R_1$ and $R_0$ as a function of $a$.
 The horizontal line represents the
value  of the Robin function $R_S$ for a round sphere of area 1.
According to \cite{okitorus} (Appendix): $R_0(a)=-\frac{\log (2 \pi)}{2\pi} - 
\frac{\log (|\eta(i a^2)|^4 a^2)}{4\pi}$ and $R_S=-\frac{1+\log\pi}{4\pi}$,
where $\eta$ is the Dedekind eta function. 
}
\label{fig2}
\end{figure}

Since $H^\prime(I)=2\pi/T(I)>0$ and $T^\prime(I)>0$ \cite{chicone1987monotonicity} (p. 315), we conclude  that
$H^{\prime\prime}(I)=-2\pi T^\prime(I)/T^2(I)<0$. This fact and  $H(0)=0$ imply
that $R_1-R_0<0$. In Figure \ref{fig2}
we present a numerical estimate
of the difference $R_1-R_0$.

The torus $\big\{\R^2/(a\Z\times a^{-1}\Z), g_1\big\}$
can be represented as the quotient of a  cylinder
that is infinite along the $x$-axis and
periodic with period $a$. We will show that this   cylinder
can be isometrically
embedded in the Euclidean three-space.
Let $X,Y,Z$ be Cartesian coordinates in $\R^3$. We will look for an
embedding of the form $X=X(x)$, $Y=F(x)\sin(2 \pi a y)$ and
$Z=F(x)\cos(2 \pi a y)$, where  $x\in\R\,,y\in \R/ a^{-1}\Z$.
The pull-back of the Euclidean metric by the embedding is
$(\dot X^2+\dot F^2)dx^2+ 4\pi^2 a^2 F^2 dy^2$. We impose that the pull-back 
 coincides with $g_1=\erm^{f(x)}(dx^2+dy^2)$ and obtain that
$4\pi^2 a^2 F^2=\erm^{f}$ and $\dot X^2+\dot F^2=\erm^{f}$. This implies
that $F(x)=\erm^{f(x)/2}/(2\pi a)$ and
\begin{equation}
  \dot X^2=\erm^{f}\left(1-\frac{\dot f^2}{16 \pi^2 a^2}\right)
  \label{X2}
\end{equation}
Since $X:\R\to\R$ must be a diffeomorphism, the right-hand side of equation
(\ref{X2}) must be strictly positive for all $x\in[0,a]$. We will show
this in the following paragraph.

The first integral (\ref{Hf}) and 
$(\erm^f-f-1)\ge 0$ imply that $\dot f^2(x)\le 2 E$ for $x\in[0,a]$,
where $E$ is the value
of $H$ for the solution with period $T(E)=a$. This  and inequality
(\ref{Taux}) imply
\[
  1-\frac{\dot f^2}{16 \pi^2 a^2}\ge 1-\frac{2 E}{16 \pi^2 T(E)^2}>0\,.
 \]

\begin{figure}[hptb!]
\centering
\begin{minipage}{0.5\textwidth}
\centering
 \includegraphics[width=0.9\textwidth]{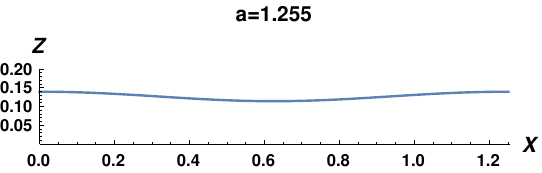}
  \includegraphics[width=0.9\textwidth]{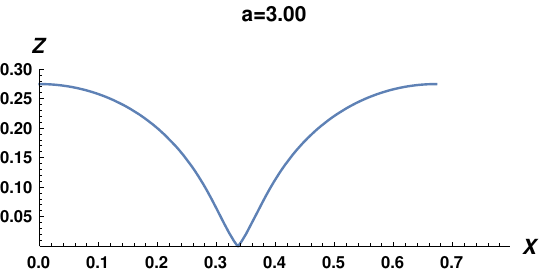}
\end{minipage}\hfill
\begin{minipage}{0.5\textwidth}
\centering
 \includegraphics[width=0.9\textwidth]{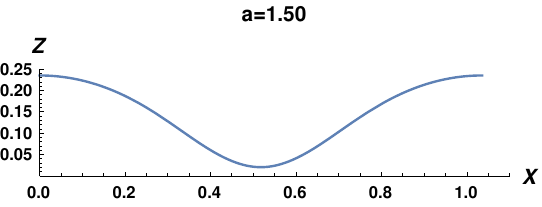}
 \includegraphics[width=0.9\textwidth]{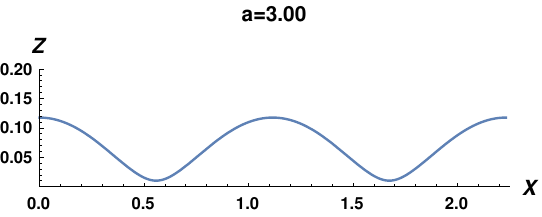}
\end{minipage}\hfill
\caption{Generating functions  of four periodic cylinders (each cylinder is
  constructed rotating the graph of  $X\to Z(X)$
  about the $X$-axis). The quotient of a cylinder
  by the group
  of periodic translations gives a torus that is isometric
  to a non-flat torus  with
  a steady-vortex metric.
   The value of the period $a$ of each torus is shown in the corresponding figure.
  There are two different tori with $a=3$: one for which the minimal period of $f$
 is $3$ and another for which the minimal period of $f$ is
 $1.5$, and  so $f$ oscillates twice inside a fundamental cell.}
\label{fig1a}
\end{figure}

\begin{figure}[hptb!]
\centering
\begin{minipage}{0.5\textwidth}
\centering
 \includegraphics[width=0.9\textwidth]{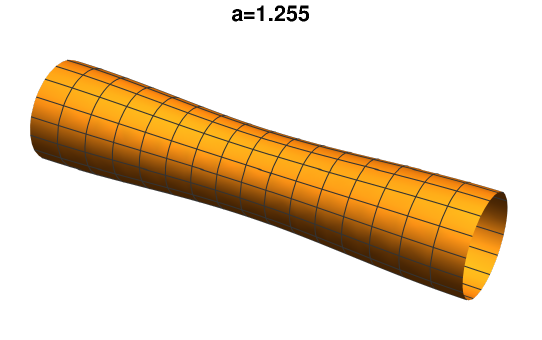}
 \includegraphics[width=0.9\textwidth]{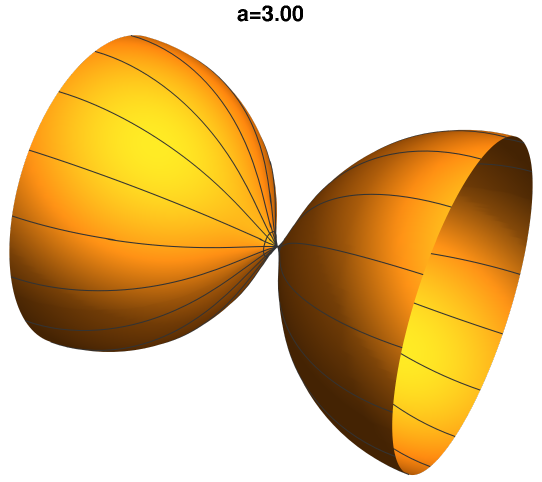}
\end{minipage}\hfill
\begin{minipage}{0.5\textwidth}
\centering
\includegraphics[width=0.9\textwidth]{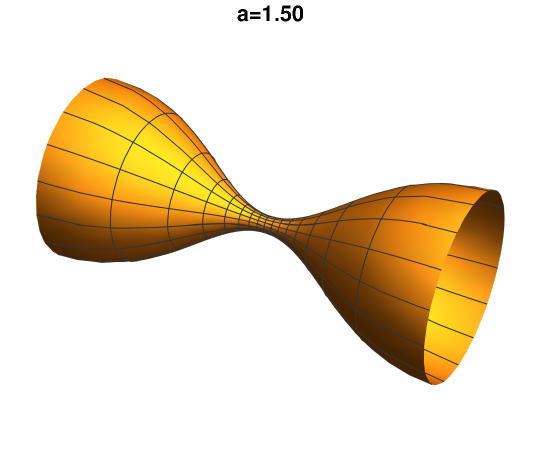}
\includegraphics[width=\textwidth]{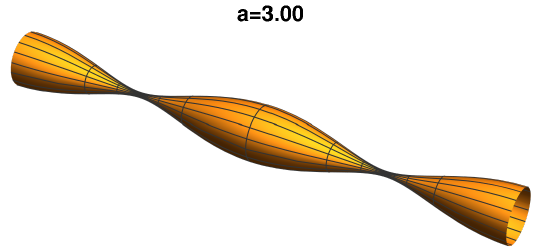}
\end{minipage}\hfill
\caption{Three-dimensional representation
  of the tori  whose generators are shown in Figure \ref{fig1a}. See
  the caption of  Figure \ref{fig1a} for explanations.}
\label{fig1b}
\end{figure}

 In  Figure \ref{fig1a} we show the  curves $x\to \{X(x), Z(x)\}$, 
 $x\in[0,a]$ and $y=0$, that when rotated about the $X-$axis generate the embedded cylinders.
 These curves were obtained by the numerical integration
 of equations  $\ddot f=8\pi(1-\erm^f)$ and
 (\ref{X2}) for: $a=1.255$, $a=1.50$, and $a=3.0$.
Only one fundamental cell of the periodic cylinder is shown.
There are two different tori with $a=3$: one for which the minimal period of $f$
 is $3$ and another for which the minimal period of $f$ is
 $1.5$, and  so $f$ oscillates twice inside a fundamental cell.
 In Figure  \ref{fig1b} we show a 3-dimensional representation of a single
 cell of each one of the cylinders whose generators are  in Figure \ref{fig1a}. 
 It is clear from  Figure \ref{fig1b} that for  $a\gg 1$ the cylinder
 becomes a collection  of aligned spheres  each one touching its  neighbors
  at a single point.
 This is in agreement with the interpretation given in
\cite{okitorus}: (the non-flat torus)  ``is approximately spherical except for a
short wormhole joining the poles''. Note:
as shown in
the right panel of Figure \ref{fig2}, in the limit as $a\to\infty$
the tori converge  to a
punctured sphere and $R_1(a)\to R_S$  where $R_S$ is the Robin function of
the round sphere.
 
\appendix

\section{Proofs of theorems \ref{crc} and \ref{sg}.}

\label{proofsec}

\begin{lemma}
\label{usmooth}
 Let $g_0$ and $g_1$ be two 
different Riemannian  metrics on $S$ in the same conformal class,
 $g_1=\lb^2g_0$ . Let $G_j$, $R_j$, $\mu_j$, $K_j$, $\Dt_j$, $j=0,1$,
  be the: 
Green's function, Robin function, volume form, 
Gaussian curvature, and Laplace operator, of $g_j$.
Let the conformal factor $\lb$ be normalized such 
that the volumes $\int_S\mu_0=\int_S\mu_1=V$ are the same. Let
$\phi$ be the unique solution of  
\[
d\ast d\phi=\frac{\mu_1-\mu_0}{V}\quad\text{with}\quad \int_S\phi\mu_0=0,
\]
that is given by 
\begin{equation}
\phi(p)=-\frac{1}{V}\int_SG_0(q,p)\lb^2(q)\mu_0(q)=-\frac{1}{V}\int_SG_0(q,p)\mu_1(q)\label{phi}
\end{equation}

Then $G_0$, $G_1$, $R_0$ and $R_1$ satisfy the following 
relations:
\begin{eqnarray}
&&G_1(q,p)-G_0(q,p)=\phi(q)+\phi(p)+c\label{t3}\\
&&\nonumber\\
&&R_1(p)=R_0(p)+\frac{1}{2\pi}\log \lambda (p)+2\phi(p)+c\label{t4}
\end{eqnarray}
where
\[
  c=
-\frac{1}{V}\int_S\phi(p)\mu_1(p)=
\frac{1}{V^2}\int_S\int_SG_0(q,p)\mu_1(q)\mu_1(p)
\]
is a constant.
Equation (\ref{t3}) is  in \cite{morpurgo1996zeta} (equation (8)) and
Equation (\ref{t4}) is in \cite{steiner2005geometrical} (Theorem 4).
\end{lemma}
\bproof
Let  $p$ and $q$ be  sufficiently close to be in 
a  domain $U$  of a 
local uniformizer $z$. 
Suppose that $U$  is  
such that any two points in $U$ are connected by a single 
geodesic in $U$. In this coordinates the length elements 
of the metrics $g_0$ and $g_1$ are $\lb_0|dz|$ and $\lb_1|dz|$, 
respectively. Notice that $\lb_1=\lb\lb_0$.
If  $\mu=dx\wedge dy$ and 
$\Dt=\frac{\partial^2}{\partial x^2}+\frac{\partial^2}{\partial y^2}$
denote the area form and   the usual Laplacian
in the coordinates
$z=(x,y)$,  respectively, then 
\begin{equation}
\Dt_j=\frac{1}{\lambda_j^2}\Dt,\qquad -\Dt\log\lambda_j=\lambda_j^2K_j,
\qquad \mu_j=\lambda^2_j \mu,\qquad j=0,1.  
\label{confr}
\end{equation}
The Dirac-delta distributions associated to the volume forms $\mu_0$ and $\mu_1$
satisfy
\[
  \delta_{j,w}=\frac{1}{\lambda^2_j}\delta_w\,,\quad\text{where}\quad
  \psi(w)=\int \psi(x,y) \delta_w(x,y)dx\wedge dy\,. 
\]

To simplify the notation we write $z(q)=z$ and 
$z(p)=w$. In the coordinates $(z,w)$ equation
(\ref{compacteq}) becomes
\begin{equation}
  -\Dt_z G_j(z,w)=\dt_w(z)-
  \frac{\lambda_j^2(z)}{V}\,.
  \label{apaux}
  \end{equation}
The Green's function can be written as
\begin{equation}
G_j(z,w)=-\frac{1}{2\pi}\log|z-w|+f_j(z,w)\label{Gzw}
\end{equation}
where $f_j(z,w)=f_j(w,z)$. 
Since $\Dt_z \log|z-w|=2\pi \dt_w$, we obtain
\begin{equation}
\Dt_z f_j(z,w)=\frac{\lb^2_j(z)}{V}. \label{Lapf}
\end{equation}

Let $\ell_j(z,w)$ be  the length with respect to the metric $g_j$
of the unique 
geodesic connecting $z$ to $w$. It can be shown that (see for
instance \cite{ragvil} proof of Theorem 5.1):
\[
\ell_j(z,w)= |w-z|\sqrt{\lb_j(z)\lb_j(w)}[1+\Oc(|z-w|)]
\]
Therefore
\[
G_j(z,w)+\frac{1}{2\pi}\log \ell_j(z,w)=
f_j(z,w)+\frac{1}{4\pi}\log[\lb_j(z)\lb_j(w)]+\Oc(|z-w|).
\]
Taking the limit as $|z-w|\to 0$ we obtain
\begin{equation}
R_j(z)=
f_j(z,z)
+\frac{1}{2\pi}\log\lb_j(z).\label{Rjz}
\end{equation}

If we subtract equation (\ref{apaux}) for $j=0$ from that for $j=1$ we obtain
\begin{equation}
 \Dt_z G_1(z,w)-\Dt_z G_0(z,w)=
  \frac{\lambda_1^2(z)-\lambda_0^2(z)}{V}\,.
  \label{apaux2}
  \end{equation}  
 This equation can be written intrinsically in terms of two-forms as  
\[
d_q\ast d_qG_1(q,p)-d_q\ast d_qG_2(q,p)=\frac{\mu_1-\mu_0}{V}
=d_q\ast d_q\phi(q)
\]
where $\phi$ is the function  in the statement of the theorem.
Thus $G_1(q,p)-G_0(q,p)=\phi(q)+ \psi(p)$ that, due to the symmetry
$G_j(p,q)=G_j(q,p)$, implies equation (\ref{t3}). Equation 
$d_q\ast d_q\phi(q)=\frac{\mu_1-\mu_0}{V}$ can be written as
$\Delta_0 \phi=(\lb^2-1)/V$. The representation formula (\ref{G})
for $\phi$ plus the relations
$\int_S\phi\mu_0=0$ and $\int_SG_0(q,p)\mu_0(q)=0$ imply that $\phi$ can be
written as in equation (\ref{phi}).
Integrating both sides of equation (\ref{t3})  with respect to
$\mu_1(q)$  over $S$ we obtain the expression for
$c$  in the lemma.
In the $z$-coordinates, 
equation (\ref{t3})  implies 
$f_1(z,w)-f_0(z,w)=\phi(z)+\phi(w)+c$. This equation  and 
equation  (\ref{Rjz}) imply 
equation (\ref{t4}).

\eproof
 \begin{lemma}
\label{conginv}
 Let $g_0$ and $g_1$ be two 
 different Riemannian metrics on $S$ in the same conformal class,
 as in Lemma {\bf \ref{usmooth}}.
Let $z=x+iy$ be a local  uniformizer and  to simplify the notation
write $z(q)=w$ and 
$z(p)=z$. Then
\begin{equation}
\left(\Dt_1R_1+\frac{K_1}{2\pi}-\frac{2}{V}\right)\mu_1=
\left(\Dt_0R_0+\frac{K_0}{2\pi}-\frac{2}{V}\right)\mu_0=-\tilde\sg\,,
\label{t5}
\end{equation}
where 
\[
-\tilde \sg= 8h(z)dx\wedge dy=4ih(z)dz\wedge d\ov z,
\]
with
\[
h(z)=  \frac{\partial}{\partial \ov w}\frac{\partial}{\partial z}
G_0(z,w)\Bigr|_{w=z}=
\frac{\partial}{\partial \ov w}\frac{\partial}{\partial z}
G_1(z,w)\Bigr|_{w=z}\,. 
\]
\end{lemma}
\bproof
In this proof we follow the notation of the proof
of Lemma \ref{usmooth}. In the $z$-coordinates,
equation (\ref{t4}) becomes
\[
R_1(z)=R_0(z)+\frac{1}{2\pi}\log \lambda_1 (z)
-\frac{1}{2\pi}\log \lambda_0 (z)+2\phi(z)+c.
\]
Taking the Laplacian $\Dt_z$ of both sides of  this equation,
 using $\Dt_z\phi_z=(\lb^2_1-\lb_0^2)/V$, and the 
 relations (\ref{confr})  for conformal metrics  we obtain 
 the first equality in equation (\ref{t5}). We recall that
$\frac{\partial}{\partial z}=
 \frac{1}{2}\left(\frac{\partial}{\partial x}
 -i\frac{\partial}{\partial y}\right)$,  
$\frac{\partial}{\partial \ov z}=
\frac{1}{2}
\left(\frac{\partial}{\partial x}+i\frac{\partial}{\partial y}\right)$,
$\Dt_z=4\frac{\partial}{\partial
  z}\frac{\partial}{\partial \ov z}$, and
 $dx\wedge dy=\frac{i}{2}dz\wedge d\ov z$.
 From equation (\ref{Rjz}) we obtain for $j=0,1$
 \begin{eqnarray*}
   \frac{\partial}{\partial \ov z}\frac{\partial}{\partial  z}
    R_j(z)&=&
\frac{\partial}{\partial \ov z}\frac{\partial}{\partial  z}
  f_j(z,w)\Bigr|_{w=z}
  +\frac{\partial}{\partial \ov w}\frac{\partial}{\partial w}
 f_j(z,w)\Bigr|_{w=z}\\
 &&+\frac{\partial}{\partial \ov z}\frac{\partial}{\partial w}
 f_j(z,w)\Bigr|_{w=z}
 +\frac{\partial}{\partial \ov w}\frac{\partial}{\partial z}
 f_j(z,w)\Bigr|_{w=z}
 +\frac{1}{2\pi}\frac{\partial}{\partial \ov z}\frac{\partial}{\partial z}
 \log\lb_j(z).
 \end{eqnarray*}
From equation (\ref{Lapf}) and $f_j(z,w)=f_j(w,z)$ we get 
\[
\frac{\partial}{\partial \ov z}\frac{\partial}{\partial  z}
f_j(z,w)\Bigr|_{w=z}=
\frac{\partial}{\partial \ov w}\frac{\partial}{\partial w}
f_j(z,w)\Bigr|_{w=z}=
 \frac{1}{4}\Dt_z f_j(z,w)\Bigr|_{w=z}=\frac{1}{4}\frac{\lb^2_j(z)}{V}
 \]
 From equation
 (\ref{Gzw})
 and from the symmetry
 $\frac{\partial}{\partial \ov z}\frac{\partial}{\partial w}
 f_j(z,w)=\frac{\partial}{\partial \ov z}\frac{\partial}{\partial w}
 f_j(w,z)$
 we get
 \begin{equation}
 \frac{\partial}{\partial \ov w}\frac{\partial}{\partial  z}
 G_j(z,w)\Bigr|_{w=z}=
\frac{\partial}{\partial \ov w}\frac{\partial}{\partial z}
f_j(z,w)\Bigr|_{w=z}
=\frac{\partial}{\partial \ov z}\frac{\partial}{\partial w}
f_j(z,w)\Bigr|_{w=z}.\label{reg}
\end{equation}
Finally,  from the above equations and from  equation (\ref{confr}) 
we obtain
\[
\Dt_z R_j(z)+\frac{\lb_j^2(z)}{2\pi}K_j(z)-\frac{2\lb^2_j(z)}{V}=
8\frac{\partial}{\partial \ov w}\frac{\partial}{\partial z}
 G_j(z,w)\Bigr|_{w=z}
\]
If we multiply both sides of this equation by $dx\wedge dy$
we obtain
\[
\left(\Dt_jR_j+\frac{K_j}{2\pi}-\frac{2}{V}\right)\mu_j=
8 \frac{\partial}{\partial \ov w}\frac{\partial}{\partial z}
G_j(z,w)\Bigr|_{w=z}dx\wedge dy
\]
for $j=0,1$. Since we have already shown that the left hand side
of this equation gives the same 2-form for $j=0$ and $j=1$, then
the right hand side has the same property.
\eproof

 The expression \(\frac{\partial}{\partial \bar{w}}\frac{\partial}{\partial z}
G_j(z,w)\) is formally analogous to the traditional Bergman kernel
for bounded domains in the complex plane. Indeed, equation
(\ref{Gzw}), which represents the decomposition of the Green's function into its singular and regular parts,
applies as well to the Green's function for bounded domains in the plane. The distinction between
the two situations lies in the regular part \( f \), which is harmonic in bounded domains,
whereas, in
this paper, the non-harmonicity of \( f \) stems from the additional term of
constant ``background vorticity.''

Following \cite{royden}, let $\partial $ be an operator
defined on complex valued functions by
$\partial=\frac{1}{2}(d+i\ast d)$ and
$\ov\partial=\frac{1}{2}(d-i\ast d)$.  In terms of a local uniformizer
$z$ we have $\partial f=\frac{\partial f}{\partial z}dz$ and
$\ov \partial f=\frac{\partial f}{\partial \ov z}d\ov z$.
\begin{lemma}
\label{bergl}
  If $G(q,p)$ is the Green's function associated to a given
  metric  and $\{\theta_1,\ldots,\theta_{2\calG}\}$ is an orthonormal basis
  of harmonic forms
  then
  \begin{equation}
  -2(\partial_p\ov\partial_qG+\ov\partial_p\partial_qG)=
  -(d_pd_qG+\ast_p\ast_qd_pd_qG)=
  \sum_{k=1}^{2\calG}\theta_k(q)\theta_k(p)
  \label{berg}
  \end{equation}
  is the Bergman reproducing kernel for harmonic forms in $S$.
  Moreover,
  if $q$ and $p$
  are in the domain of a local uniformizer with $z(q)=w$ and $z(p)=z$
  then
  \begin{equation}
  2(\partial_p\ov\partial_qG+\ov\partial_p\partial_qG)
  =4\Real \{\partial_p\ov\partial_qG\}=
4\Real\left\{\frac{\partial}{\partial z}\frac{\partial}{\partial \ov w}
 G(w,z)d\ov w dz\right\}
  \label{bergl2}
  \end{equation}
  
\end{lemma}
\bproof
The equality
$2(\partial_p\ov\partial_qG+\ov\partial_p\partial_qG)=
d_pd_qG+\ast_p\ast_qd_pd_qG$ and  equation (\ref{bergl2})
are direct consequences of the definition
of the operators $\partial$ and $\ov\partial$. 
Due to equation (\ref{reg}) the function
$\frac{\partial}{\partial z}\frac{\partial}{\partial \ov w}
G(w,z)$ is $\cinfty$ for all values of $z$ and $w$ including $z=w$.
So the double one-form $d_pd_qG+\ast_p\ast_qd_pd_qG$ is $\cinfty$
on $S\times S$.

The Bergman reproducing kernel for harmonic forms
$H(q,p)=\sum_{k=1}^{2\calG}\theta_k(q)\theta_k(p)$ is characterized by the
following properties:
\begin{itemize}
\item[]  For an arbitrary function $\psi$ on $S$:
  \begin{eqnarray*}
  && \quad\int_{S}d\psi(p)\wedge H(q,p)=0\\
  && \quad \int_{S}\ast_pd\psi(p)\wedge H(q,p)=0 \,, \\
  \end{eqnarray*}
  where the integrations are with respect to  the variable $p$;
  and for any harmonic one-form $\nu$ on $S$
\[
\nu(q)=\int_{S}\nu(p)\wedge \ast_p H(q,p)=
 \sum_{k=1}^{2\calG}\theta_k(q)\int_{S(p)}\nu(p)\wedge\ast_p\theta_k(p).
\]
\end{itemize}

In order to  prove the equality $d_pd_qG+\ast_p\ast_qd_pd_qG=-H(q,p)$
we use the regularity of $d_pd_qG+\ast_p\ast_qd_pd_qG$ on $S\times S$.
So, for any function $\psi$ on $S$ 
\[
\begin{split}
 & \int_{S}d_p\psi(p)\wedge(d_pd_qG+\ast_p\ast_qd_pd_qG)=
  -\int_{S}\psi(p)\wedge d_p(d_pd_qG+\ast_p\ast_qd_pd_qG)\\
  &\qquad\qquad =-\lim_{\ep\to 0}\int_{S-B_\ep(q)}\psi(p)\wedge
  d_p(d_pd_qG+\ast_p\ast_qd_pd_qG)
  \end{split}
\]
where $B_\ep(q)$ is a small ball (with respect to any local uniformizer)
of radius $\ep$ with center at $q$. For $p$ outside $B_\ep(q)$,  
\[
d_p(d_pd_qG+\ast_p\ast_qd_pd_qG)=\ast_qd_q(d_p\ast_pd_p G)=
\ast_qd_q\left(\frac{\mu(p)}{V}\right)=0,
\]
so $\int_{S}d_p\psi(p)\wedge(d_pd_qG+\ast_p\ast_qd_pd_qG)=0$.
In the same way it is possible to prove that
$\int_{S}\ast_pd\psi(p)\wedge(d_pd_qG+\ast_p\ast_qd_pd_qG)=0$.

It remains to show that
$
\nu(q)=-\int_{S}\nu(p)\wedge \ast_p (d_pd_qG+\ast_p\ast_qd_pd_qG)
$
for any harmonic one-form $\nu$ on $S$.
This is a consequence of
\[
\begin{split}
&  
\int_{S}\nu(p)\wedge \ast_p (d_pd_qG+\ast_p\ast_qd_pd_qG)
=\lim_{\ep\to 0}\int_{S-B_\ep(q)}\nu(p)\wedge
  \ast_p(d_pd_qG+\ast_p\ast_qd_pd_qG)\\
&\qquad =\lim_{\ep\to 0}
   \int_{-\partial B_\ep(q)}\ast_p\nu(p)d_qG+
  \int_{-\partial B_\ep(q)}\nu(p)\ast_qd_qG
\end{split}
  \]
  An explicit computation using a local uniformizer gives that
  this last integral is equal to $-\nu(q)$.
  \eproof


Theorem \ref{sg} is a consequence of lemmas \ref{conginv} and 
\ref{bergl} and the following reasoning.
Let $z(p)=z=x+iy$ and $z(q)=w=\xi+i\eta$
be the components of the local uniformizer used in lemma
\ref{bergl} and $\theta_k(p)=\theta_{k1}(z)dx+\theta_{k2}(z)dy$
and $\theta_k(q)=\theta_{k1}(w)d\xi+\theta_{k2}(w)d\eta$
be the components of $\theta_k$. Lemma \ref{bergl} implies that
\[
\begin{split}
  \sum_{k=1}^{2\calG}\theta_k(q)\theta_k(p)=&
\left(\sum_{k=1}^{2\calG}\theta_{k1}(w)\theta_{k1}(z)\right)dxd\xi+
\left(\sum_{k=1}^{2\calG}\theta_{k2}(w)\theta_{k2}(z)\right)dyd\eta\\
&
+\left(\sum_{k=1}^{2\calG}\theta_{k2}(w)\theta_{k1}(z)\right)dxd\eta+
\left(\sum_{k=1}^{2\calG}\theta_{k1}(w)\theta_{k2}(z)\right)dyd\xi\\
&=-4\Real\left\{\frac{\partial}{\partial z}
\frac{\partial}{\partial \ov w} G(w,z)d\ov w dz\right\}
\end{split}
\]
For $q=p$ and  $dz=dw$,  the right hand side of this  equation becomes
\[
-4
\frac{\partial}{\partial z}
\frac{\partial}{\partial \ov w} G(w,z)\Bigr|_{w=z}(dx^2+dy^2)
\]
that implies
\[
\sum_{k=1}^{2\calG}\theta^2_{k1}(z)=
\sum_{k=1}^{2\calG}\theta^2_{k2}(z)=
-4
\frac{\partial}{\partial z}
\frac{\partial}{\partial \ov w} G(w,z)\Bigr|_{w=z},\quad\text{and}\quad
\sum_{k=1}^{2\calG}\theta_{k1}(z)\theta_{k2}(z)=0
\]
So, the  form $\sg$ in theorem \ref{sg} can be written as
\[
\begin{split}
\sg(z)&=\sum_{k=1}^{2\calG} \theta_k(z)\wedge *\theta_{k}(z)=
\sum_{k=1}^{2\calG} [\theta_{k1}(z)dx+\theta_{k2}(z)dy]\wedge
    [\theta_{k1}(z)dy-\theta_{k2}(z)dx]\\
&=\sum_{k=1}^{2\calG} [\theta^2_{k1}(z)+\theta^2_{k2}(z)]dx\wedge dy
 =-8
\frac{\partial}{\partial z}
\frac{\partial}{\partial \ov w} G(w,z)\Bigr|_{w=z}dx\wedge dy=\tilde \sg,   
\end{split}
\]
where $\tilde \sg$ is the form in equation (\ref{t5}).
This proves that equation (\ref{Req}) holds and
finishes the proof of theorem \ref{sg}.
\eproof

Now we prove theorem \ref{crc}.
The Robin function on a Riemannian manifold $(S,g)$ is constant 
whenever $(S,g)$     admits a transitive Lie group action  of 
isometries. So, the Robin function is constant for the round sphere
 and for all flat tori.
Let $S$ be a   sphere (torus) endowed with a  Riemannian metric  $g_0$.
The  uniformization  theorem  implies the
existence of  a diffeomorphism from 
$(S,g_0)$ to the round sphere (a flat torus) $(\S^2,g_1)$ 
 such that the pull-back of $g_1$ is conformal to $g_0$. So, the 
existence of a steady vortex metric on the sphere (torus) is
proved.

The proof is  more complicated 
when $S$ is compact and has a genus larger than one. 
Equation (\ref{t4}) implies:
\begin{equation}
\Dt_0 R_1(p)=\Dt_0R_0(p)+\frac{1}{2\pi}\Dt_0\log \lambda (p)+2
\frac{\lb^2-1}{V}
\label{eqeq}
\end{equation}
Imposing that $R_1$ is constant, normalizing the volume $V$  of $S$ 
to be equal to one, and defining 
\[
u=4\pi R_0+\log\lb^2
\]
we get the  following equation for $u$
\begin{equation}
\Dt_0u=8\pi-8\pi h\erm^u \label{higgs}
\end{equation}
where $h=\erm^{-4\pi R_0}$.
To each solution of this equation corresponds  a Riemannian 
metric $g_1$  conformal to $g_0$ such that $\Dt R_1=0$ and therefore
$R_1$  is constant.
Equation (\ref{higgs}) was very much studied for several 
reasons. It appears in the problem of finding a Riemannian 
metric on the sphere with a prescribed curvature $h$ 
that is conformal to the standard metric 
with curvature $4\pi$ (the conformal factor is $\erm^u$).
It also appears  in the so-called Chern-Simons-Higgs theory
(see \cite{dingjost} for references). The following theorem was 
taken from \cite{dingjost} (it is 
a combination of their theorem 1.2 plus their
remark 1.3).
\begin{theorem}[Ding, Jost, Li,and  Wang]
\label{ding}
Let $(S,g_0)$ be a compact Riemann surface and  let $K_0$ be its
Gauss curvature. Let $h$ be a positive smooth function on $S$.
 Suppose that the function 
$8\pi R_0 + 2 \log h$  achieves its maximum at $p$.
 If $\Dt_0 \log h(p) > -(8\pi - 2K_0(p))$ then 
equation {\rm (\ref{higgs})}  has a smooth solution.
\end{theorem}
It is remarkable that in the case we are interested in  
$h=\erm^{-4\pi R_0}$ and $8\pi R_0 + 2 \log h=0$. So, any point  
in $S$ is a point of maximum and therefore to finish the 
proof it is sufficient 
to show the existence of  a point $p$ in $S$ where the inequality 
$0>-\Dt_0 \log h(p)-(8\pi - 2K_0(p))$ holds.
The Gauss-Bonet theorem implies 
$\int_SK_0\mu_0=2\pi (2-2\calG)$, where $\calG$ is the genus of $S$.
Since $\int \mu_0=1$, 
 the integral of the right hand side of the inequality above is
$-8\pi \calG<0$. This finishes the proof of existence of a natural 
vortex metric if $\calG>1$. 
\eproof

\section{The Robin function and
  the Minakshisundaram–Pleijel zeta function.}
\label{min}

The Minakshisundaram–Pleijel zeta function, which will be referred as
the zeta function, is defined as
\begin{equation}\label{zetapqs}
  \zeta(q,p,s)=\sum_{k=1}^\infty\frac{\phi_k(q)\phi_k(p)}{\lambda_k^s}=
  \frac{1}{\Gamma(s)}\int_0^\infty\left(K(q,p,t)-\frac{1}{V}\right)t^{s-1}dt\,,
 \end{equation}
 where  $s\in\C$ and  Re $s>n/2$
 (the convergence is a consequence of  inequality (\ref{asymp1})).

 According to the theorem in Section 5 of  \cite{MP},
 the function
 \[
   \zeta(p,s):=\zeta(p,p,s)
 \]
 can be extended as a meromorphic function to the whole complex plane. If
 dimension $n\ge 3$
is odd, then the only possible poles of $\zeta(p,s)$ are located at
$s=n/2,n/2-1,\ldots,3/2,1/2,-1/2,\ldots$.
If the dimension $n$ is even,
then $\zeta(p,s)$ has at most a finite number of poles
   that are possibly located at   
   $s=n/2,n/2-1,\ldots,2,1$ and 
 the residue at the poles can be computed \cite{MP}. In particular, if $n$ is even and $s$ is close to $s=1$, then  
 \begin{equation}
   \zeta(p,s)=\frac{1}{(4\pi)^{n/2}}\frac{a_{n/2-1}(p)}{s-1} +\text{convergent power
   series in }(s-1)\,,\label{pole1}
\end{equation}
where $a_{n/2-1}(p)$ is the function that appears in equation (\ref{asymp2}).

If $s$ is made equal to one in equation (\ref{zetapqs}), then  we obtain a
formal expression
\begin{equation}\label{formal}
  G(q,p)= \int_0^\infty\left(K(q,p,t)-\frac{1}{V}\right)dt=
  \sum_{k=1}^\infty\frac{\phi_k(q)\phi_k(p)}{\lambda_k} \ {}^\prime\!\!=^\prime\zeta(q,p,1)\,
\end{equation}
that indicates a possible relation between the regularization of
$G(q,p)$ and $\zeta(q,p,s)$ as $q\to p$ and $s\to 1$.
Indeed, for $n=2$ the following result holds (see, e.g.
\cite{steiner2005geometrical}, Proposition 2 and the Appendix):
\begin{equation}\begin{split}
 R(p)=&
  \lim_{\ell(q,p)\to 0} \left[
    G(q,p) +\frac{1}{2\pi}\log \ell(q,p) \right]
  \\ =&\lim_{s\to 1} \left[\zeta(p,s)-\frac{1}{(4\pi)}\frac{1}{s-1}\right]
  +\frac{\log 4-2 \gamma}{4\pi}\end{split}\label{steiner}
\end{equation}
where
$\gamma$ is the Euler's constant. In the following theorem we show that this result
can be generalized to higher dimensions.
The same  result, for an  elliptic operator
that appears in the context of quantum field theory in curved spacetime, 
was obtained by Bilal and Ferrari in \cite{bilal2013multi} (Section 3).
If the parameters $m$ and $\psi$ that appear in their elliptic operator
are set equal to zero, then the formulas  in  equations (3.45) and (3.46) of
\cite{bilal2013multi} are exactly ours in theorem (\ref{Robinzeta}). 

\begin{theorem}\label{Robinzeta}
  The Robin function can be written in terms of the analytic extension of the 
  Minakshisundaram–Pleijel zeta function as
  \begin{equation}
    \begin{array}{rcll}
   R(p)&=&\lim_{s\to 1} \left[\zeta(p,s)-\frac{1}{(4\pi)}\frac{1}{s-1}\right]
  +\frac{\log 4-2 \gamma}{(4\pi)^{n/2}} & \text{if $n$ is even,}\\ & & & \\
      R(p)&=&\zeta(p,1) &\text{if $n$ is odd}\,.
    \end{array}
    \end{equation}
\end{theorem}
\bproof We will  prove the theorem only  for  $n$ even, since the proof
for $n$ odd is similar. For $s>n/2$ both sides of equation (\ref{zetapqs})
converge. The idea is to add terms to both sides of that equation
such that the integral in the right-hand side of  equation (\ref{zetapqs})
converges when $s=1$. In analogy to what we did to define the Robin function
we  rewrite equation (\ref{zetapqs}) for $s>n/2$ as
\begin{equation}\label{zetapps}\begin{split}
    & \zeta(p,s)-\sum_{k=0}^{n/2-1}\frac{a_k(p)}{(4\pi)^{n/2}\Gamma(s)}\int_0^1t^{k-n/2}t^{s-1}dt =\frac{1}{\Gamma(s)}\int_1^\infty\left(K(p,p,t)-\frac{1}{V}
    \right)t^{s-1}dt\\
 &\quad +\lim_{\epsilon\to 0_+}
 \frac{1}{\Gamma(s)}\int_\epsilon^1\left(K(p,p,t)-\frac{1}{V}
-\sum_{k=0}^{n/2-1}\frac{a_k(p)}{(4\pi)^{n/2}}t^{k-n/2}
 \right)t^{s-1}dt\,.\end{split}
 \end{equation}
 For $s>n/2$, the left-hand side of this equation can be written as
 \begin{equation}
   \zeta(p,s)-\frac{a_{n/2-1}(p)}{(4\pi)^{n/2}\Gamma(s)}\frac{1}{s-1}-
   \sum_{k=0}^{n/2-2}\frac{a_k(p)}{(4\pi)^{n/2}\Gamma(s)}\frac{1}{s-n/2+k}\,.
 \end{equation}
 Due to equations (\ref{asymp1}) and (\ref{asymp2}), the integrand in the last
 line of equation  (\ref{zetapps}) is bounded by a constant times $t^{s-1}$, and
 therefore the right-hand side of equation (\ref{zetapps}) is an analytic
 function of $s$ for Re$\, s>0$. This implies that the analytic continuation
 of $\zeta(p,s)$ to  Re$\, s>0$ is given by the regular function at the right-hand
 side of equation (\ref{zetapps}) plus the poles at $s=1,2,\ldots,n/2$
 explicitly given in the left-hand side of the same equation. With this
 understanding,  we can compute the regularized value of $\zeta(p,s)$ at $s=1$ as
 \begin{equation}\label{intaux2}
   \begin{split}
     &  \lim_{s\to 1} \left[ \zeta(p,s)-\frac{a_{n/2-1}(p)}{(4\pi)^{n/2}\Gamma(s)}
       \frac{1}{s-1}\right]\\ & =
     \lim_{s\to 1}\left[ \zeta(p,s)-\frac{a_{n/2-1}(p)}{(4\pi)^{n/2}}
       \frac{1}{s-1}\right]+\frac{a_{n/2-1}(p)}{(4\pi)^{n/2}}\Gamma^\prime(1)\\
     &= \sum_{k=0}^{n/2-2}\frac{a_k(p)}{(4\pi)^{n/2}}\frac{1}{1-n/2+k}
     +\int_1^\infty\left(K(p,p,t)-\frac{1}{V}
    \right)dt\\
 &\quad +\lim_{\epsilon\to 0_+}
 \int_\epsilon^1\left(K(p,p,t)-\frac{1}{V}
-\sum_{k=0}^{n/2-1}\frac{a_k(p)}{(4\pi)^{n/2}}t^{k-n/2}
 \right)dt\,,\end{split}
     \end{equation}
   where we  used  that  the integrand in the last
   line of equation  (\ref{zetapps}) is bounded by a constant times $t^{s-1}$ to
   exchange the order of the limits. Performing the integrals of the terms
   that are polynomials in $t$ in the right-hand side of equation (\ref{intaux2}),
   using the definition of  the Robin function given in equation (\ref{Re}),
   and that $\Gamma^\prime(1)=-\gamma$ we
   obtain the result in the statement of the theorem.
   \eproof

\nd {\bf Acknowledgments.}
This paper is dedicated to  Jair Koiller who introduced me to the subject
of vortices on surfaces and presented to me the work of Okikiolu and Steiner.
Jair  has been a constant source of inspiration.

\vskip 1cm

\nd
Data sharing not applicable to this article as no datasets were generated or analysed during the current study.
\bibliographystyle{plain}

\end{document}